\documentclass[preprint,12pt]{elsarticle}

\usepackage{amstext,amssymb,amsmath,amsbsy,bm}
\usepackage{hyperref}
\usepackage{amscd}
\usepackage{amsfonts}
\usepackage{indentfirst}
\usepackage{verbatim}
\usepackage{amsmath}
\usepackage{amsthm}
\usepackage{enumerate}
\usepackage{graphicx}
\usepackage{color}
\usepackage[OT1]{fontenc}
\usepackage[latin1]{inputenc}
\usepackage[english]{babel}
\usepackage{amssymb}
\usepackage{subfig}
\usepackage{algorithm}
\usepackage{algpseudocode}
\usepackage[shortlabels]{enumitem}

\newcommand{\x}{{\bf x}}

\newcommand{\bmu}{\boldsymbol{\mu}}

\newcommand{\bu}{{\bf u}}
\newcommand{\bU}{{\bf U}}
\newcommand{\bF}{{\bf F}}
\newcommand{\bv}{{\bf v}}
\newcommand{\bw}{{\bf w}}
\newcommand{\bW}{{\bf W}}

\newcommand{\Div}{{\rm div}}

\newtheorem{Theorem}{Theorem}

\newtheorem{Proposition}{Proposition}
\newtheorem{Corollary}{Corollary}
\newtheorem{Remark}{Remark}

\newtheorem*{Assumption*}{Assumption}

\newtheorem{Problem}{Problem}
\newtheorem*{Problem*}{Problem}
\numberwithin{equation}{section}

\usepackage{matlab-prettifier}
\usepackage{mathtools}

\journal{}

\begin{document}

\begin{frontmatter}

\title{The inverse initial data problem for anisotropic Navier--Stokes equations via Legendre time reduction method} 

\author[label1]{Cong B. Van} 
\ead{cvan1@charlotte.edu}
\cortext[cor1]{Corresponding author.}
\affiliation[label1]{organization={Department of Mathematics and Statistics, University of North Carolina at
Charlotte},
            city={Charlotte},
            postcode={28223}, 
            state={NC},
            country={USA}}

\author[label2]{Thuy T. Le} 
\ead{Thuy.Le@csulb.edu}

\affiliation[label2]{organization={Department of Mathematics and Statistics, 
California State University, Long Beach},
            city={Long Beach},
            postcode={90840}, 
            state={CA},
            country={USA}}
            
\author[label1]{Loc H. Nguyen\corref{cor1}} 
\ead{loc.nguyen@charlotte.edu.}

\begin{abstract}
We consider an inverse initial-data problem for the compressible anisotropic Navier--Stokes equations, in which the goal is to reconstruct the initial velocity field from noisy lateral boundary observations. In the formulation studied here, the density, pressure, anisotropic viscosity tensor, and body force are assumed known, while the initial velocity is the quantity to be recovered. We introduce a new computational framework based on Legendre time-dimensional reduction, in which the velocity field is projected onto an exponentially weighted Legendre basis in time. This transformation reduces the original time-dependent inverse problem to a coupled system of time-independent elliptic equations for the Fourier coefficients of the velocity field. The resulting reduced model is solved using a combination of quasi-reversibility and a damped Picard iteration. Numerical experiments in two dimensions show that the proposed method accurately and robustly reconstructs initial velocity fields, even in the presence of significant measurement noise, geometrically complex structures, and anisotropic effects. The method provides a flexible and computationally tractable approach for inverse fluid problems in anisotropic media.
\end{abstract}




\end{frontmatter}



\section{Introduction}

Let $d \geq 2$ denote the spatial dimension, and let $\Omega \subset \mathbb{R}^d$ be a bounded domain with smooth boundary $\partial \Omega$. We consider a fourth-order tensor field $\boldsymbol{\mu}: \overline{\Omega} \to \mathbb{R}^{d \times d \times d \times d}$ representing the viscosity property of the medium,
  which satisfies the coercivity condition: there exists a constant $\Lambda > 0$ such that for all matrices $\boldsymbol{\xi} = (\xi_{ij})_{i, j = 1}^d \in \mathbb{R}^{d \times d}$, we have
    \begin{equation*}
        \sum_{i,j,k,l=1}^d \mu_{ijkl} \, \xi_{ij} \, \xi_{kl} \geq \Lambda |\boldsymbol{\xi}|^2.
    \end{equation*}
Let $T>0$ be a final time. We consider a fluid velocity field
$\bu:\Omega\times(0,T)\to\mathbb{R}^d$ governed by the following
initial--boundary value problem for the compressible anisotropic Navier--Stokes system:
\begin{equation}
    \begin{cases}
        \displaystyle
        \rho \left( \frac{\partial \bu}{\partial t}
        + (\bu \cdot \nabla)\bu \right)
        = -\nabla p + \Div \left( \bmu : \nabla \bu \right)
        + \rho(\x, t)\bF(\x,t),
        & (\x,t) \in \Omega \times (0,T),\\[0.4em]
        \rho_t + \Div (\rho \bu) = 0,
        & (\x,t) \in \Omega \times (0,T),\\[0.4em]
        p = c^2 \rho
        \qquad \text{or} \qquad
        p = K\rho^\gamma,
        & (\x,t) \in \Omega \times (0,T),\\[0.4em]
        \bu(\x,t) = 0,
        & (\x,t) \in \partial\Omega \times (0,T),\\[0.4em]
        \bu(\x,0) = \bu^0(\x),
        & \x \in \Omega,\\[0.4em]
        \rho(\x,0)=\rho^0(\x),
        & \x \in \Omega.
    \end{cases}
    \label{main_eqn}
\end{equation}Here, $\bu^0:\Omega\to\mathbb{R}^d$ denotes the initial velocity field,
$\rho:\Omega\times(0,T)\to\mathbb{R}$ is the density, $\rho^0(\x)$ is a given positive initial density,
$p:\Omega\times(0,T)\to\mathbb{R}$ is the pressure,
and $\bF:\Omega\times(0,T)\to\mathbb{R}^d$ is the body force.
The constants $c>0$, $K>0$, and $\gamma>1$ are given, where
$p=c^2\rho$ corresponds to the isothermal case and
$p=K\rho^\gamma$ corresponds to the polytropic law. 
 
Our focus is on the following inverse problem:
\begin{Problem}[Inverse initial data problem]\label{isp}
Given the lateral boundary observation
\begin{equation}
\mathbf{f}(\x, t) = \partial_{\nu}\mathbf{u}(\x, t)
\quad \mbox{for all } (\x, t) \in \partial \Omega \times (0, T),
\label{data}
\end{equation}
assume that the density $\rho(\x,t)$, the pressure $p(\x,t)$, the anisotropic viscosity tensor $\bmu$,
and the body force $\bF(\x,t)$ are known.
Reconstruct the initial condition $\bu^0(\x)$ in $\Omega$.
\end{Problem}

To avoid ambiguity in the formulation of the inverse problem, we now clarify which quantities are intended to be reconstructed and which are treated as known quantities. Since the observation \eqref{data} provides only boundary information, attempting to reconstruct multiple independent space-time fields from this data would be severely underdetermined without supplementary measurements or additional structural assumptions. We emphasize that the simultaneous recovery of quantities such as $\rho(\x,t)$, $p(\x,t)$, and $\bF(\x,t)$ would be mathematically significant and practically very interesting. However, in the present setting, the observation is given only on a $(d-1)$-dimensional spatial manifold over time, making the total dimension $(d-1)+1=d$, whereas these unknown fields are defined in the full $d$-dimensional space-time cylinder $\Omega\times(0,T)$, whose dimension is $d+1$. Therefore, the available data do not contain enough information to identify all these quantities simultaneously. Accordingly, we formulate a targeted inverse problem in which the initial condition $\bu^0$ is the sole reconstruction target, while $\rho(\x,t)$, $p(\x,t)$, $\bmu$, and $\bF(\x,t)$ are assumed known.

The inverse initial data problem for anisotropic Navier--Stokes equations holds substantial practical importance in a variety of scientific and engineering domains, including seismology, aerospace and automotive engineering, hydraulic systems, oceanography, and atmospheric modeling~\cite{Lions1998, Stokes1851, Tsai2018, White1991}. In these contexts, direct access to the initial velocity field is often infeasible, while accurate reconstruction of this field is essential for initializing predictive simulations of fluid dynamics. Recovering the initial state from lateral boundary measurements thus constitutes a critical and timely inverse problem.
From a theoretical standpoint, this problem presents significant challenges. The well-posedness of the forward anisotropic Navier--Stokes system remains unresolved in general settings, particularly in three dimensions and without strong regularity assumptions. For the inverse problem, the absence of Carleman estimates or other analytical tools applicable to fully anisotropic fourth-order tensors makes it difficult to establish uniqueness or stability. As a result, many foundational theoretical questions remain open in this setting.

In this work, we focus exclusively on the computational aspect of the problem. We do not attempt to address the underlying analytical difficulties but instead develop a numerical framework based on time-dimensional reduction and iterative regularization. The proposed method achieves stable and accurate reconstructions in practice, providing a viable computational approach for a class of inverse problems that currently lack rigorous theoretical guarantees.

While the primary focus of this paper is on an inverse problem associated with the anisotropic Navier--Stokes equations, it is essential first to understand the status of the corresponding forward problem. The well-posedness and regularity of the forward model form the foundation upon which the theoretical analysis and numerical treatment of inverse problems are built.  Therefore, we briefly review relevant results on the forward problem.
While the global existence and regularity of solutions to the Navier--Stokes equations remain one of the most challenging open problems in mathematical fluid dynamics, substantial progress has been made in specific settings. Foundational analytical frameworks have been developed for various formulations of the Navier--Stokes system, particularly in the compressible and stationary cases \cite{Feireisl2004,Ladyzhenskaya1969, Lions1998,NovotnyStraskraba2004, Temam1979}. The existence of global weak solutions for stationary problems with homogeneous Dirichlet boundary conditions was first established in~\cite{Lions1998}, and later refined and generalized in~\cite{NovotnyStraskraba2004,FrehseGojSteinhauer2005,PlotnikovSokolowski2004,PlotnikovSokolowski2005}. 
Further developments in the well-posedness theory include the treatment of compressible stationary Navier--Stokes equations with small data, where uniqueness and differentiability were proven in~\cite{PlotnikovRubanSokolowski2008}. More recent advances have addressed coupled systems and complex boundary interactions: global well-posedness for a Navier--Stokes--Cahn--Hilliard model with chemotaxis and singular potential in two dimensions was established in~\cite{HeWu2021}, while local and global results for the Boussinesq system with viscoelastic effects under various boundary conditions were obtained in~\cite{AntontsevKhompysh2023}. In three dimensions, global existence and decay properties of strong solutions to the nonhomogeneous incompressible Navier--Stokes--Cahn--Hilliard system were demonstrated in~\cite{FangNeiGuo2025} using a priori energy methods. Finally, the well-posedness of stationary Navier--Stokes equations in higher-dimensional half-spaces was studied in~\cite{Fujii2025}, where scaling-critical Besov spaces and maximal regularity theory were employed to derive existence and asymptotic behavior results.

Inverse problems for the Navier--Stokes equations have been extensively studied to determine unknown coefficients or initial data from boundary or interior measurements. In \cite{FouresteyMoubachir2005}, a numerical approximation for inverse problems was developed using a Lagrange--Galerkin method that effectively addresses drag reduction and velocity identification. A regularized Gauss--Newton method for shape reconstruction in two-dimensional Navier--Stokes equations was derived in \cite{YanHeMa2010}, establishing differentiability with respect to the boundary curve. Global identifiability of viscosity in two-dimensional Stokes and Navier--Stokes equations using Cauchy force measurements was proved in \cite{LaiUhlmannWang2015}, advancing prior Dirichlet-to-Neumann map approaches. Global uniqueness for inverse boundary value problems for the Navier--Stokes equations and Lam\'e system in two dimensions was established in \cite{ImanuvilovYamamoto2015a} using complex geometric optics solutions. Iterative methods, including a novel simple backward integration approach, were introduced in \cite{OConnorEtAl2024} to solve retrospective inverse problems, outperforming direct adjoint looping. Global uniqueness for viscosity determination in Stokes and Navier--Stokes equations in two and three dimensions was demonstrated in \cite{Liu2024}, combining prior theorems with new results for arbitrary bounded domains.

Unlike the standard PDE-constrained optimization framework based on state, adjoint, and optimality systems, our approach follows a direct reconstruction route. Our strategy for solving the inverse problem associated with the Navier--Stokes equations stated in Problem~\ref{isp} is grounded in the use of the Legendre time-reduction technique. Specifically, we approximate the velocity field $\mathbf{u}$ in the time variable by truncating its Fourier series with respect to an exponentially weighted Legendre basis of the weighted space $L_{e^{-2t}}^2(0,T)$. This basis, recently introduced in~\cite{TrongElastic} for inverse initial data problems in elasticity, and further utilized in~\cite{LeVanDangNguyen} for Maxwell's system, offers notable advantages over classical bases such as the trigonometric system and the Legendre polynomial basis without the exponential weight. In particular, the exponentially weighted Legendre basis possesses two desirable properties: nonvanishing derivatives and favorable convergence behavior of the truncated series applied to time derivatives of the velocity field. By applying this time-reduction procedure, we reformulate the original time-dependent inverse problem into a coupled system of time-independent equations, where the principal operator takes the form $\Div(\boldsymbol{\mu}:\nabla\cdot)$. We refer to this reformulated system as the time-dimensional reduction model. Leveraging this structure, we then employ quasi-reversibility together with a damped Picard iteration to solve the reduced model. The obtained solution directly yields an approximate solution to the original inverse problem, thereby providing a practical iterative reconstruction scheme.

The remainder of the paper is organized as follows. 
Section~\ref{sec:legendre_review} reviews the exponentially weighted Legendre basis used for temporal approximation. Section~\ref{sec:reduction} derives the time-dimensional reduction model that converts the time-dependent inverse problem into a coupled elliptic system. Section~\ref{sec:Picard} presents a Picard, least-squares (quasi-reversibility) algorithm for solving the reduced model. Section~\ref{sec:numerics} reports numerical experiments demonstrating accuracy and robustness. We then include a dedicated Section \ref{sec:discuss} on the stability of the time reduction and the choice of time basis, including a scaled weight for large final times $T$. Section~\ref{sec:conclusion} concludes with a discussion and future directions.

\section{The Legendre polynomial-exponential basis}
\label{sec:legendre_review}

The Legendre polynomial-exponential basis, first introduced in \cite{TrongElastic}, is fundamental to the time-dimensional reduction method employed in this study. This basis integrates the spectral characteristics of classical Legendre polynomials with an exponential weight, providing a structured approximation framework in the space of exponentially weighted square-integrable functions.  
For clarity and completeness, we briefly revisit this basis below.
Let $\{P_n\}_{n \geq 0}$ denote the Legendre polynomials on $(-1,1)$, defined by Rodrigues' formula:
\[
P_n(x) = \frac{1}{2^n n!} \frac{d^n}{dx^n}(x^2 - 1)^n.
\]
To define an orthonormal basis on $(0, T)$, we introduce the affine transformation $x = \frac{2t}{T} - 1$, and define the rescaled polynomials
\[
Q_n(t) := \sqrt{\frac{2n + 1}{T}} P_n\left(\frac{2t}{T} - 1\right) \quad \mbox{for } t \in (0, T).
\]
The set $\{Q_n\}_{n \geq 0}$ forms an orthonormal basis in $L^2(0, T)$.
We then define the Legendre polynomial-exponential basis by:
\[
\Psi_n(t) := e^t Q_n(t) \quad \mbox{for } t \in (0, T).
\]
The system $\{\Psi_n\}_{n \geq 0}$ forms an orthonormal basis of the exponentially weighted space
\[
L^2_{e^{-2t}}(0, T) := \left\{ u \in L^2(0, T) \,\middle|\, \int_0^T e^{-2t} |u(t)|^2 dt < \infty \right\},
\]
with the inner product
\[
\langle u, v \rangle_{L^2_{e^{-2t}}(0, T)} := \int_0^T e^{-2t} u(t) v(t) \, dt.
\]

\begin{Proposition}
The Legendre polynomial-exponential basis functions $\Psi_n$, for $n \geq 0$, satisfy the following regularity and spectral properties:
\begin{enumerate}
    \item (Smoothness and non-vanishing derivatives) For each $n \geq 0$, the function $\Psi_n$ is infinitely differentiable on $(0, T)$, and none of its derivatives of any order vanishes identically on this interval.
    
    \item (Spectral coefficient decay) For every integer $k \in \mathbb{N}$, there exists a constant $C > 0$, depending only on $k$ and $T$, such that for all $u \in H^k(0, T)$,
    \begin{equation} \label{Hk-inequality}
        \sum_{n = 0}^\infty n^{2k} \left| \langle u, \Psi_n \rangle_{L^2_{e^{-2t}}(0, T)} \right|^2 \leq C \|u\|_{H^k(0, T)}^2.
    \end{equation}
    
    \item (Derivative norm growth) There exists a constant $C > 0$, depending only on $T$, such that for all $n \geq 1$,
    \begin{equation}
        \|\Psi_n'\|_{L^2_{e^{-2t}}(0, T)} \leq C n^{3/2},
        \quad \text{and} \quad
        \|\Psi_n''\|_{L^2_{e^{-2t}}(0, T)} \leq C n^{7/2}.
        \label{3.2}
    \end{equation}
\end{enumerate}
\label{pro2_1}
\end{Proposition}

\begin{Remark}[The importance of the weight $e^t$]
We omit the proof of Proposition~\ref{pro2_1}, which follows from Proposition~2.1, Lemma~2.1, and Lemma~2.2 in \cite{TrongElastic}. 
Note that Item~1 relies on the factor $e^{t}$ in $\Psi_n=e^{t}Q_n$: without this factor, derivatives of some $\Psi_n$ may be identically zero (e.g., the constant time mode), which can reduce the effectiveness of the time-dimensional reduction; see Subsection~\ref{subsec:exponent} for details.
\end{Remark}

Let $N$ be a positive integer. We introduce three key operators associated with this basis.

\begin{itemize}
    \item \textbf{Fourier Coefficient Operator} ($\mathbb{F}^N$): For  $u \in L^2_{e^{-2t}}\left((0, T); L^2(\Omega)^d\right)$, we define the vector of its first $N+1$ coefficients by:
    \[
    \mathbb{F}^N[u](\x) := \begin{bmatrix} u_0(\x) & u_1(\x) & \cdots & u_N(\x) \end{bmatrix}^\top,
    \]
    where
    \[
    u_n(\x) := \int_0^T e^{-2t} u(\x, t) \Psi_n(t) dt
    \]
    for $\x \in \Omega.$
    
    \item \textbf{Expansion Operator} ($\mathbb{S}^N$): Given the coefficient vector 
    \[         
        \mathbf{U} = (u_0, \ldots, u_N)^\top \in  L^2(\Omega)^{d \times (N + 1)},     
    \]
    the space-time expansion is defined by
    \[
    \mathbb{S}^N[\mathbf{U}](\x, t) := \sum_{n=0}^N u_n(\x) \Psi_n(t),
    \]
    for $(\x, t) \in \Omega \times (0, T).$
    \item \textbf{Projection Operator} ($\mathbb{P}^N$): The projection of $u \in L^2_{e^{-2t}}((0, T); L^2(\Omega)^{d})$ onto the subspace $\mathbb{V}^N := \mathrm{span}\{\Psi_0, \ldots, \Psi_N\}$ is:
    \begin{equation}
    \mathbb{P}^N[u](\x, t) := \sum_{n=0}^N \left( \int_0^T e^{-2t} u(\x, t) \Psi_n(t) dt \right) \Psi_n(t),
    \label{P}
    \end{equation}
    for $(\x, t) \in \Omega \times (0, T).$
\end{itemize}

These operators satisfy the identity
\begin{equation*}
\mathbb{S}^N[\mathbb{F}^N[u]] = \mathbb{P}^N[u],
\end{equation*}
and Parseval's identity implies:
\begin{equation}
\|\mathbb{P}^N[u]\|_{L^2_{e^{-2t}}((0, T); L^2(\Omega)^d)} = \|\mathbf{U}\|_{L^2(\Omega)^{d \times (N+1)}} \leq \|u\|_{L^2_{e^{-2t}}((0, T); L^2(\Omega)^d)}.
\label{Par}
\end{equation}
\begin{Remark}
    For any integer $N \geq 1$, the operators $\mathbb{F}^N$, $\mathbb{S}^N$, and $\mathbb{P}^N$ admit natural restrictions to vector-valued functions in the space $L^2_{e^{-2t}}\left((0, T); H^p(\Omega)^d\right)$ for any $p \geq 1$. In this setting, the Parseval-type inequality in \eqref{Par} continues to hold with $L^2(\Omega)^d$ replaced by $H^p(\Omega)^d$. Explicitly, we have
    \begin{equation*}
        \|\mathbb{P}^N[u]\|_{L^2_{e^{-2t}}((0, T); H^p(\Omega)^d)} = \|\mathbf{U}\|_{H^p(\Omega)^{d \times (N+1)}} \leq \|u\|_{L^2_{e^{-2t}}((0, T); H^p(\Omega)^d)}.        
    \end{equation*}
\end{Remark}

\section{The time-dimensional reduction method} \label{sec:reduction}
The inverse problem stated in Problem~\ref{isp} can be reformulated as the recovery of a vector field $\bu$ satisfying
\begin{equation}
\begin{cases}
    \displaystyle
    \rho(\x,t)\left(\partial_t \bu + (\bu \cdot \nabla)\bu \right)
    = - \nabla p(\x,t) + \Div (\bmu : \nabla \bu) + \rho(\x,t)\bF(\x,t),
    & \text{in } \Omega \times (0, T), \\[0.4em]
    \bu(\x, t) = \mathbf{0},
    & \text{on } \partial \Omega \times (0, T), \\[0.4em]
    \partial_{\nu} \bu(\x, t) = \mathbf{f}(\x, t),
    & \text{on } \partial \Omega \times (0, T),
\end{cases}
\label{4_1}
\end{equation}
where $\rho(\x,t)$, $p(\x,t)$, $\bmu$, and $\bF(\x,t)$ are known, while the initial condition $\bu^0$ is the unknown to be reconstructed.

\begin{Remark}
The existence of a solution to \eqref{4_1} is justified by the fact that $\bu$
is the restriction of the solution to the forward problem \eqref{main_eqn}
on $\Omega \times (0, T)$. A rigorous analysis of uniqueness for the original
inverse problem is beyond the scope of this paper, whose main focus is on the
computational methodology. Nevertheless, this question remains of independent
theoretical interest.
\end{Remark}

\begin{Remark}
The continuity equation and any closure relation between $p$ and $\rho$
(such as $p=c^2\rho$ or $p=K\rho^\gamma$) belong to the forward compressible model.
However, they are not explicitly needed in the reconstruction step considered in this paper.
Indeed, in our targeted inverse problem, the density $\rho(\x,t)$ and the pressure $p(\x,t)$
are assumed known, together with $\bmu$ and $\bF(\x,t)$. Therefore, the momentum equation alone
forms a closed system of $d$ equations for the $d$ components of the unknown velocity field $\bu$.
For this reason, the time-dimensional reduction procedure is applied only to the momentum equation.
\end{Remark}

We now derive the time-dimensional reduction model. Let $N \geq 1$ be a fixed integer, and consider the projection operator $\mathbb{P}^N$. Applying $\mathbb{P}^N$ to both sides of the differential equation in \eqref{4_1}, we obtain
\begin{equation}
    \mathbb{P}^N\!\left[\rho(\x,t)\partial_t \bu\right]
    +
    \mathbb{P}^N\!\left[\rho(\x,t)(\bu \cdot \nabla)\bu\right]
    =
    - \mathbb{P}^N[\nabla p(\x,t)]
    + \Div \left( \bmu : \nabla \mathbb{P}^N[\bu] \right)
    + \mathbb{P}^N\!\left[\rho(\x,t)\bF(\x,t)\right]
    \label{4_3333}
\end{equation}
in $\Omega \times (0, T)$,
where we used the fact that $\bmu$ is known and independent of $t$.
 In \eqref{4_3333}, we have used the identity
\[
\mathbb{P}^N[\Div(\bmu:\nabla \bu)]
=
\Div\bigl(\bmu:\nabla \mathbb{P}^N[\bu]\bigr).
\]
This is true because $\mathbb{P}^N$ acts only on the temporal variable and $\Div(\bmu:\nabla,\cdot,)$ is a spatial operator with time-independent coefficient tensor $\bmu$.

Adding and subtracting
$\rho(\x,t)\partial_t\mathbb{P}^N[\bu]
+\rho(\x,t)(\mathbb{P}^N[\bu]\cdot\nabla)\mathbb{P}^N[\bu]$
on the left-hand side of \eqref{4_3333}, and then rearranging the resulting equation, we obtain
\begin{multline}
    \rho(\x,t)\partial_t \mathbb{P}^N[\bu]
    - \Div \left( \bmu : \nabla \mathbb{P}^N[\bu] \right)
    + \rho(\x,t)(\mathbb{P}^N[\bu] \cdot \nabla)\mathbb{P}^N[\bu]
    \\
    + \mathbb{P}^N[\nabla p(\x,t)]
    - \mathbb{P}^N[\rho(\x,t)\bF(\x,t)]
    =
    -\Bigl(\mathbb{P}^N[\rho(\x,t)\partial_t \bu]
    - \rho(\x,t)\partial_t \mathbb{P}^N[\bu]\Bigr)
    \\
    -\Bigl(\mathbb{P}^N[\rho(\x,t)(\bu \cdot \nabla)\bu]
    - \rho(\x,t)(\mathbb{P}^N[\bu] \cdot \nabla)\mathbb{P}^N[\bu]\Bigr)
    \label{4_3}
\end{multline}
in $\Omega \times (0, T)$.

We have the theorem
\begin{Theorem}[Asymptotic commutation theorem]\label{thm1}
Assume that the solution $\bu$ to \eqref{4_1} satisfies
$
\bu \in H^3\bigl((0,T);H^s(\Omega)^d\bigr)
$
for some integer $s>\frac d2+1$. Then
\begin{equation}
\lim_{N\to\infty}
\left\|
\mathbb{P}^N[\partial_t \bu]-\partial_t\mathbb{P}^N[\bu]
\right\|_{L^2((0,T);H^s(\Omega)^d)}
=0,
\label{3,4}
\end{equation}
and
\begin{equation}
\lim_{N\to\infty}
\left\|
\mathbb{P}^N[(\bu\cdot\nabla)\bu]
-(\mathbb{P}^N[\bu]\cdot\nabla)\mathbb{P}^N[\bu]
\right\|_{L^2((0,T);H^{s-1}(\Omega)^d)}
=0.
\label{3,5}
\end{equation}
\end{Theorem}

\begin{proof}
Since the weighted norm on $L^2_{e^{-2t}}(0,T)$ is equivalent to the standard
$L^2(0,T)$ norm on the finite interval $(0,T)$, we will work with the standard
time norms below.

We first prove \eqref{3,4}. By the approximation property of the projection
operator $\mathbb{P}^N$ associated with the weighted Legendre basis, we have
\[
\mathbb{P}^N[\bu]\to \bu
\quad\text{in }H^1\bigl((0,T);H^s(\Omega)^d\bigr),
\]
and
\[
\mathbb{P}^N[\partial_t\bu]\to \partial_t\bu
\quad\text{in }L^2\bigl((0,T);H^s(\Omega)^d\bigr).
\]
Therefore,
\begin{align*}
\left\|
\mathbb{P}^N[\partial_t \bu]-\partial_t\mathbb{P}^N[\bu]
\right\|_{L^2((0,T);H^s(\Omega)^d)}
&\le
\left\|
\mathbb{P}^N[\partial_t \bu]-\partial_t\bu
\right\|_{L^2((0,T);H^s(\Omega)^d)}
\\
&\quad
+
\left\|
\partial_t\bu-\partial_t\mathbb{P}^N[\bu]
\right\|_{L^2((0,T);H^s(\Omega)^d)}.
\end{align*}
Both terms on the right-hand side converge to zero as $N\to\infty$, which proves
\eqref{3,4}.

We next prove \eqref{3,5}. Define the nonlinear mapping
\[
\mathcal N(\bv):=(\bv\cdot\nabla)\bv.
\]
Since $s>\frac d2$, the product of two functions in $H^s(\Omega)$ still belongs to $H^s(\Omega)$.
Moreover, since $s>\frac d2+1$, the mapping
\[
\mathcal N:H^s(\Omega)^d\to H^{s-1}(\Omega)^d
\]
is locally Lipschitz. In particular, there exists a constant $C>0$ such that for all
$\bv,\bw\in H^s(\Omega)^d$,
\begin{equation}
\|\mathcal N(\bv)-\mathcal N(\bw)\|_{H^{s-1}(\Omega)^d}
\le
C\bigl(\|\bv\|_{H^s(\Omega)^d}+\|\bw\|_{H^s(\Omega)^d}\bigr)
\|\bv-\bw\|_{H^s(\Omega)^d}.
\label{eq:lipschitz_nonlinear}
\end{equation}
Now write
\begin{align*}
&\left\|
\mathbb{P}^N[(\bu\cdot\nabla)\bu]
-(\mathbb{P}^N[\bu]\cdot\nabla)\mathbb{P}^N[\bu]
\right\|_{L^2((0,T);H^{s-1}(\Omega)^d)}
\\
&\qquad\le
\left\|
\mathbb{P}^N[(\bu\cdot\nabla)\bu]
-(\bu\cdot\nabla)\bu
\right\|_{L^2((0,T);H^{s-1}(\Omega)^d)}
\\
&\qquad\quad+
\left\|
(\bu\cdot\nabla)\bu
-(\mathbb{P}^N[\bu]\cdot\nabla)\mathbb{P}^N[\bu]
\right\|_{L^2((0,T);H^{s-1}(\Omega)^d)}.
\end{align*}
The first term tends to zero because
\[
(\bu\cdot\nabla)\bu \in L^2\bigl((0,T);H^{s-1}(\Omega)^d\bigr)
\]
and $\mathbb{P}^N$ converges to the identity in this space.

For the second term, applying \eqref{eq:lipschitz_nonlinear} pointwise in time and then
integrating over $(0,T)$, we obtain
\begin{align*}
&\left\|
(\bu\cdot\nabla)\bu
-(\mathbb{P}^N[\bu]\cdot\nabla)\mathbb{P}^N[\bu]
\right\|_{L^2((0,T);H^{s-1}(\Omega)^d)}
\\
&\qquad\le
C\Bigl(
\|\bu\|_{L^\infty((0,T);H^s(\Omega)^d)}
+
\|\mathbb{P}^N[\bu]\|_{L^\infty((0,T);H^s(\Omega)^d)}
\Bigr)
\|\bu-\mathbb{P}^N[\bu]\|_{L^2((0,T);H^s(\Omega)^d)}.
\end{align*}
Since
\[
\mathbb{P}^N[\bu]\to \bu
\quad\text{in }H^1\bigl((0,T);H^s(\Omega)^d\bigr),
\]
the Sobolev embedding
\[
H^1\bigl((0,T);H^s(\Omega)^d\bigr)\hookrightarrow
C\bigl([0,T];H^s(\Omega)^d\bigr)
\]
implies that $\{\mathbb{P}^N[\bu]\}_{N\ge1}$ is uniformly bounded in
$L^\infty((0,T);H^s(\Omega)^d)$, and
\[
\|\bu-\mathbb{P}^N[\bu]\|_{L^2((0,T);H^s(\Omega)^d)}\to 0.
\]
Hence, the second term also converges to zero. This proves \eqref{3,5}.
\end{proof}

\begin{Remark}[Asymptotic commutation]
For large $N$, due to Theorem \ref{thm1}, the spectral projection $\mathbb{P}^N$ approximatedly commutes with $\partial_t$ and with $(\bu\cdot\nabla)\bu$. In other words, ``project-then-operate" and ``operate-then-project" agree in the limit. This consideration motivates our choice of the exponentially weighted Legendre basis $\{\Psi_n\}_{n\ge 0}$. Other time bases (e.g., the exponentially weighted polynomial basis in $L^2(0,T)$; see \cite{Klibanov:jiip2017,NguyenLeNguyenKlibanov:2023}) need not have this property and may therefore be unsuitable for deriving \eqref{4.2}.
\end{Remark}

\begin{Corollary}\label{cor:weighted_commutation}
Assume that $\rho \in C^1(\overline{\Omega}\times[0,T])$ is known and that
\[
\bu \in H^3\bigl((0,T);H^s(\Omega)^d\bigr)
\]
for some integer $s>\frac d2+1$. Then
\begin{equation}
\lim_{N\to\infty}
\left\|
\mathbb{P}^N[\rho(\x,t)\partial_t \bu]
-
\rho(\x,t)\partial_t\mathbb{P}^N[\bu]
\right\|_{L^2((0,T);H^s(\Omega)^d)}
=0,
\label{cor_rho_dt}
\end{equation}
and
\begin{equation}
\lim_{N\to\infty}
\left\|
\mathbb{P}^N[\rho(\x,t)(\bu\cdot\nabla)\bu]
-
\rho(\x,t)(\mathbb{P}^N[\bu]\cdot\nabla)\mathbb{P}^N[\bu]
\right\|_{L^2((0,T);H^{s-1}(\Omega)^d)}
=0.
\label{cor_rho_nl}
\end{equation}
\end{Corollary}

\begin{proof}

To prove \eqref{cor_rho_dt}, we note first that
\[
\mathbb{P}^N[\rho(\x,t)\partial_t \bu]
\to
\rho(\x,t)\partial_t \bu
\quad \text{in } L^2((0,T);H^s(\Omega)^d),
\]
by the approximation property of $\mathbb{P}^N$. On the other hand, by \eqref{3,4},
\[
\partial_t\mathbb{P}^N[\bu]\to \partial_t\bu
\quad \text{in } L^2((0,T);H^s(\Omega)^d).
\]
Since $\rho(\x,t)$ is smooth, multiplication by $\rho(\x,t)$ is continuous on
$H^s(\Omega)^d$. Therefore,
\[
\rho(\x,t)\partial_t\mathbb{P}^N[\bu]
\to
\rho(\x,t)\partial_t\bu
\quad \text{in } L^2((0,T);H^s(\Omega)^d).
\]
Combining these two convergences, we obtain \eqref{cor_rho_dt}.

The proof of \eqref{cor_rho_nl} is analogous. Indeed,
to prove \eqref{cor_rho_nl}, we first note that
\[
\mathbb{P}^N[\rho(\x,t)(\bu\cdot\nabla)\bu]
\to
\rho(\x,t)(\bu\cdot\nabla)\bu
\quad \text{in } L^2((0,T);H^{s-1}(\Omega)^d),
\]
by the approximation property of $\mathbb{P}^N$. On the other hand, by \eqref{3,5},
\[
(\mathbb{P}^N[\bu]\cdot\nabla)\mathbb{P}^N[\bu]
\to
(\bu\cdot\nabla)\bu
\quad \text{in } L^2((0,T);H^{s-1}(\Omega)^d).
\]
Again, since multiplication by the smooth coefficient $\rho(\x,t)$ is continuous,
we obtain
\[
\rho(\x,t)(\mathbb{P}^N[\bu]\cdot\nabla)\mathbb{P}^N[\bu]
\to
\rho(\x,t)(\bu\cdot\nabla)\bu
\quad \text{in } L^2((0,T);H^{s-1}(\Omega)^d).
\]
Hence \eqref{cor_rho_nl} follows.
\end{proof}

\begin{Remark}[Regularity in Theorem~\ref{thm1}, Corollary~\ref{cor:weighted_commutation}, and observed performance]
Theorem~\ref{thm1} and Corollary~\ref{cor:weighted_commutation} require relatively high regularity of the solution. In particular, we assume
\[
\bu \in H^3\bigl((0,T);H^s(\Omega)^d\bigr)
\]
for some integer $s>\frac d2+1$, and in the corollary we additionally assume that the known density $\rho(\x,t)$ is sufficiently smooth. These assumptions are introduced to justify the projection-based analysis and to obtain clear convergence statements for both the unweighted and $\rho$-weighted terms.

In computations, however, the algorithm performs well even when the initial condition $\bu^0$ belongs only to $L^2(\Omega)^d$ and may be discontinuous, so the assumptions of Theorem~\ref{thm1} and Corollary~\ref{cor:weighted_commutation} are not satisfied; see Section~\ref{sec:numerics} for our out-of-expectation numerical results. This suggests that the practical scope of the method may be substantially broader than what is currently justified by the available theory.
\end{Remark}

Due to Corollary~\ref{cor:weighted_commutation} and the limits
\eqref{cor_rho_dt}--\eqref{cor_rho_nl}, the right-hand side of \eqref{4_3}
vanishes asymptotically as $N\to\infty$. Therefore, equation~\eqref{4_3}
can be approximated by the reduced system
\begin{multline}
    \rho(\x,t)\partial_t \mathbb{P}^N[\bu]
    - \Div \left( \bmu : \nabla \mathbb{P}^N[\bu] \right)
    + \rho(\x,t)(\mathbb{P}^N[\bu] \cdot \nabla)\mathbb{P}^N[\bu]
    \\
    + \mathbb{P}^N[\nabla p(\x,t)]
    - \mathbb{P}^N[\rho(\x,t)\bF(\x,t)]
    = \mathbf{0},
    \label{4_6}
\end{multline}
for $(\x, t) \in \Omega \times (0, T),$
which serves as the foundation of our time-dimensional reduction model.

For each $n \geq 0$, let
\begin{equation}
    \bu_n(\x) = \int_0^T e^{-2t}\bu(\x,t)\Psi_n(t)\,dt,
    \label{4.2}
\end{equation}
be the $n$-th Fourier mode of $\bu$ with respect to the basis
$\{\Psi_n\}_{n\geq 0}$.

Due to \eqref{P},
\begin{equation}
   \mathbb{P}^N [\bu](\x, t)
    =
    \sum_{n = 0}^N \bu_n(\x)\Psi_n(t).
    \label{u_app}
\end{equation}
Substituting \eqref{u_app} into the system \eqref{4_6}, we obtain the approximation
\begin{multline}
\rho(\x,t)\sum_{n=0}^N \bu_n(\x)\Psi_n'(t)
-\sum_{n=0}^N \Div \left( \bmu : \nabla \bu_n(\x) \right)\Psi_n(t)
\\
+\rho(\x,t)\sum_{l=0}^N \sum_{n=0}^N
\left( \bu_l(\x) \cdot \nabla \right) \bu_n(\x)\Psi_l(t)\Psi_n(t)
+\mathbb{P}^N[\nabla p(\x,t)]
-\mathbb{P}^N[\rho(\x,t)\bF(\x,t)]
= \mathbf{0},
\label{4.3}
\end{multline}
valid for all $(\x,t)\in\Omega\times(0,T)$.

To derive the time-reduced model, we test equation \eqref{4.3} against
$e^{-2t}\Psi_m(t)$ and integrate over $t\in(0,T)$ for each $m=0,1,\dots,N$.
This yields
\begin{multline}
\sum_{n=0}^N s_{mn}(\x)\bu_n(\x)
+\sum_{l=0}^N\sum_{n=0}^N a_{mnl}(\x)
\left(\bu_l(\x)\cdot\nabla\right)\bu_n(\x)
\\
=
-\nabla p_m(\x)
+\Div\left(\bmu:\nabla \bu_m(\x)\right)
+\bF_m(\x),
\label{4.5}
\end{multline}
for all $\x\in\Omega$, where
\[
s_{mn}(\x)=\int_0^T e^{-2t}\rho(\x,t)\Psi_n'(t)\Psi_m(t)\,dt,
\]
\[
a_{mnl}(\x)=\int_0^T e^{-2t}\rho(\x,t)\Psi_l(t)\Psi_n(t)\Psi_m(t)\,dt,
\]
and
\[
p_m(\x)=\int_0^T e^{-2t}p(\x,t)\Psi_m(t)\,dt,
\qquad
\bF_m(\x)=\int_0^T e^{-2t}\rho(\x,t)\bF(\x,t)\Psi_m(t)\,dt.
\]

The corresponding Cauchy boundary conditions for each mode $\bu_m$ follow from
\eqref{data} and \eqref{4.2}:
\begin{equation}
    \bu_m(\x) = \mathbf{0},
    \qquad
    \partial_{\nu}\bu_m(\x)
    =
    \int_0^T e^{-2t}\mathbf{f}(\x,t)\Psi_m(t)\,dt,
    \label{4.6}
\end{equation}
for all $\x \in \partial \Omega$.

Equations~\eqref{4.5} and~\eqref{4.6} define a coupled boundary value problem
for the vector of coefficients $\mathbb{P}^N[\bu]$. We refer to this system as the
\emph{time-dimensional reduction model}:
\begin{equation}
    \begin{cases}
\displaystyle
\sum_{n=0}^N s_{mn}(\x)\bu_n(\x)
+\sum_{l=0}^N\sum_{n=0}^N a_{mnl}(\x)
\left(\bu_l(\x)\cdot\nabla\right)\bu_n(\x)
\\
\hspace{5cm}=
-\nabla p_m(\x)
+\Div\left(\bmu:\nabla\bu_m(\x)\right)
+\bF_m(\x),
& \x \in \Omega,\\[0.8em]
\displaystyle
\bu_m(\x)=\mathbf{0},
& \x \in \partial\Omega,\\[0.4em]
\displaystyle
\partial_{\nu}\bu_m(\x)
=
\mathbf{f}_m(\x)
:=
\int_0^T e^{-2t}\mathbf{f}(\x,t)\Psi_m(t)\,dt,
& \x \in \partial\Omega,
    \end{cases}
    \label{time_red}
\end{equation}
for $m=0,1,\dots,N$.

Solving the system~\eqref{time_red} yields the coefficient functions
$\bu_0,\bu_1,\dots,\bu_N$, which can be used to reconstruct an approximation of the original
space-time solution $\bu(\x,t)$ via
\[
\mathbb{P}^N[\bu](\x,t)=\sum_{n=0}^N \bu_n(\x)\Psi_n(t).
\]
According to the approximation formula~\eqref{u_app}, this reconstruction provides
a finite-dimensional approximation of the original space-time solution $\bu$.
Hence, the inverse problem under consideration (Problem~\ref{isp}) is transformed into the
computational task of solving the reduced model~\eqref{time_red}.

Despite this reduction, solving the system~\eqref{time_red} remains highly nontrivial due to the presence of the nonlinear convective term
\begin{equation}
\sum_{l=0}^N \sum_{n=0}^N a_{mnl}(\x)
\left( \bu_l(\x) \cdot \nabla \right)\bu_n(\x),
\qquad
m\in\{0,1,\dots,N\},
\label{convection_term}
\end{equation}
which makes \eqref{time_red} a coupled nonlinear system.
To overcome this difficulty, we apply a conventional damped Picard iteration; see Section~\ref{sec:Picard}.

\begin{Remark}
The time-dimensional reduction model~\eqref{time_red} is an approximate representation of the original system, arising from the truncation of the infinite Fourier expansion in~\eqref{u_app}. Nevertheless, we adopt the working assumption that
\[
\mathbf{U}
=
\begin{bmatrix}
\bu_0 & \bu_1 & \dots & \bu_N
\end{bmatrix}^{\top}
\]
satisfies the reduced model.
This assumption is justified by Theorem~\ref{thm1} and Corollary~\ref{cor:weighted_commutation}, which show that the projection and weighted commutation errors vanish as $N\to\infty$.
Consequently, the reduced model becomes an increasingly accurate surrogate for the full dynamics as the cutoff parameter $N$ increases.
From a computational viewpoint, this assumption is also natural, since any practical implementation necessarily uses a finite truncation level $N$.
\end{Remark}

\begin{Remark}
The uniqueness of the time-reduction model~\eqref{time_red} remains an open problem due to the complexity introduced by the anisotropic fourth-order tensor $\bmu$ and the nonlinear convection term \eqref{convection_term}. A standard approach in inverse problems is to use Carleman estimates to establish uniqueness and stability. However, deriving such estimates in the fully anisotropic setting is highly nontrivial and, to the best of our knowledge, is not yet available in a general form. In contrast, in the isotropic case, Carleman estimates are more tractable; see, for instance, \cite{ImanuvilovLorenziYamamoto2022, Isakov2007}. Under such simplifications, the uniqueness of the time-reduction model~\eqref{time_red} may become accessible to rigorous analysis. We leave this direction for future research.
\end{Remark}

\section{A damped Picard algorithm to solve the time-reduction model}\label{sec:Picard}

We define the admissible set
\[
H := \left\{
    \mathbf{U} \in H^2(\Omega)^{d \times (N + 1)}
    \ \middle| \
    \mathbf{U}|_{\partial \Omega} = \mathbf{0}
\right\},
\]
which serves as the admissible space for the reduced system in~\eqref{time_red}. Let the initial guess be
\[
{\bf U}^{(0)} =
\begin{bmatrix}
\bu^{(0)}_0 & \bu^{(0)}_1 & \dots & \bu^{(0)}_N
\end{bmatrix}^{\top}.
\]
Assume that
\[
{\bf U}^{(k)} =
\begin{bmatrix}
\bu^{(k)}_0 & \bu^{(k)}_1 & \dots & \bu^{(k)}_N
\end{bmatrix}^{\top}
\]
is known for some $k\geq 0$. We first compute an intermediate iterate
\[
\widetilde{\bf U}^{(k+1)} =
\begin{bmatrix}
\widetilde{\bu}^{(k+1)}_0 & \widetilde{\bu}^{(k+1)}_1 & \dots & \widetilde{\bu}^{(k+1)}_N
\end{bmatrix}^{\top}
\]
as the solution of the following linearized boundary value problem:
\begin{equation}
    \begin{cases}
\displaystyle
\sum_{n=0}^N s_{mn}(\x)\widetilde{\bu}_n^{(k+1)}(\x)
-\Div\left(\bmu:\nabla \widetilde{\bu}_m^{(k+1)}(\x)\right)
\\[0.4em]
\displaystyle\hspace{2cm}
=
-\sum_{l=0}^N\sum_{n=0}^N a_{mnl}(\x)
\left(\bu_l^{(k)}(\x)\cdot\nabla\right)\bu_n^{(k)}(\x)
-\nabla p_m(\x)+\bF_m(\x),
& \x\in\Omega,\\[0.8em]
\widetilde{\bu}_m^{(k+1)}(\x)=\mathbf{0},
& \x\in\partial\Omega,\\[0.4em]
\partial_{\nu}\widetilde{\bu}_m^{(k+1)}(\x)=\mathbf{f}_m(\x),
& \x\in\partial\Omega,
    \end{cases}
    \label{linear}
\end{equation}
for each $m=0,1,\dots,N$.

We then define the damped Picard update by
\begin{equation}
{\bf U}^{(k+1)}
=
(1-\omega){\bf U}^{(k)}
+
\omega\,\widetilde{\bf U}^{(k+1)},
\qquad \omega=0.5.
\label{damped_picard}
\end{equation}
This iteration is a damped fixed-point scheme, where the nonlinear convection term is evaluated at the previous iterate, and the damping parameter is used to improve stability.

The next step is to solve the over-determined problem \eqref{linear}. By over-determined, we mean that both Dirichlet and Neumann boundary conditions are prescribed for $\widetilde{\bf U}^{(k+1)}$. We propose to use the quasi-reversibility method, first introduced in \cite{LattesLions:e1969}. 
Thus, the intermediate iterate $\widetilde{\bf U}^{(k+1)}$ is defined as the minimizer of the least-squares functional
\begin{multline}
    J_{{\bf U}^{(k)}}^{\epsilon,N}(\boldsymbol{\varphi})
    =
    \sum_{m=0}^N \Bigg[
    \int_{\Omega}
    \Bigg|
    \Div\left(\bmu:\nabla \varphi_m(\x)\right)
    - \sum_{n=0}^N s_{mn}(\x)\varphi_n(\x)
\\
    - \sum_{l=0}^N \sum_{n=0}^N a_{mnl}(\x)
    \left(\bu_l^{(k)}(\x)\cdot\nabla\right)\bu_n^{(k)}(\x)
     -\nabla p_m(\x)
    + \bF_m(\x)
    \Bigg|^2 d\x
\\
    + \int_{\partial\Omega}
    \left|\partial_{\nu}\varphi_m(\x)-\mathbf{f}_m(\x)\right|^2 d\sigma(\x)
    + \epsilon \|\varphi_m\|_{H^2(\Omega)^d}^2
    \Bigg].
    \label{5.3}
\end{multline}
The existence and uniqueness of the minimizer of $J_{{\bf U}^{(k)}}^{\epsilon, N}$ follow from the strict convexity of the functional and the fact that its domain $H$ is a Hilbert space.

We summarize the damped Picard iteration and quasi-reversibility method for solving \eqref{time_red}, and hence Problem~\ref{isp}, in Algorithm~\ref{alg}.

\begin{algorithm}[h!]
\caption{\label{alg}Solution procedure for Problem~\ref{isp} via the time-dimensional reduction method}
\begin{algorithmic}[1]
\State  Select a cutoff number $N$ and the regularization parameter $\epsilon$. \label{s1}
\State Set the damping parameter $\omega = 0.5$.
\State Set ${\bf U}^{(0)}={\bf 0}\in H$.
\For{ \label{s4} $k=0$ to $K-1$, for some $K\in\mathbb{N}$}
    \State Minimize the functional $J_{{\bf U}^{(k)}}^{\epsilon,N}$ defined in \eqref{5.3}, and denote the obtained minimizer by $\widetilde{\bf U}^{(k+1)}$.  \label{s5}
    \State Update the iterate by the damped Picard step
    \[
    {\bf U}^{(k+1)}
    =
    (1-\omega){\bf U}^{(k)}
    +
    \omega\,\widetilde{\bf U}^{(k+1)}.
    \]
\EndFor \label{s7}
\State Set the computed solution ${\bf U}_{\rm comp}={\bf U}^{K}$:
\[
{\bf U}_{\rm comp}=
\begin{bmatrix}
\bu_{0, {\rm comp}} & \bu_{1, {\rm comp}} & \dots & \bu_{N, {\rm comp}}
\end{bmatrix}^{\top}.
\]
\State Reconstruct the fluid velocity field $\bu_{\rm comp}$ using the space-time expansion
\begin{equation}
\bu_{\rm comp}(\x,t)
=
\mathbb{S}^N[{\bf U}_{\rm comp}]
=
\sum_{n=0}^N \bu_{n, {\rm comp}}(\x)\Psi_n(t),
\label{4.3333}
\end{equation}
for all $(\x,t)\in\Omega\times(0,T)$.
\State Compute the approximate initial velocity field
\[
\bu^0_{\rm comp}(\x)=\bu_{\rm comp}(\x,0),
\qquad \x\in\Omega.
\]
\end{algorithmic}
\end{algorithm}

\begin{Remark}
While our numerical experiments indicate that the iterative sequence $\{{\bf U}^{(k)}\}_{k\geq 0}$ converges reliably to the target solution when initialized from zero, we do not provide a theoretical justification of this convergence. The main difficulty lies in the absence of an appropriate Carleman estimate tailored to the anisotropic Navier--Stokes equations in the time-reduced setting. Developing such an estimate is highly nontrivial due to the complex structure of the fourth-order viscosity tensor. In the isotropic case, where the viscosity tensor has a simpler form, one may hope to analyze a similar iterative scheme using existing Carleman estimates; see, for instance, \cite{ImanuvilovLorenziYamamoto2022, Isakov2007}. Since the present work focuses on anisotropic media, we leave the corresponding theoretical convergence analysis for future research.
\end{Remark}

\section{Numerical study}\label{sec:numerics}

In this section, we present numerical experiments to validate the proposed time-dimensional reduction method and to illustrate its performance in reconstructing the initial velocity field from noisy boundary observations.

\subsection{Manufactured data generation}

We present numerical simulations in the two-dimensional setting ($d=2$) for simplicity. 
In our experiments, the fourth-order viscosity tensor $\bmu$ is chosen to be spatially constant, and its flattened representation is given by
\[
\boldsymbol{\mu}_{\rm flat}
=
\begin{bmatrix}
\mu_{1111} & \mu_{1112} & \mu_{1121} & \mu_{1122}\\
\mu_{1211} & \mu_{1212} & \mu_{1221} & \mu_{1222}\\
\mu_{2111} & \mu_{2112} & \mu_{2121} & \mu_{2122}\\
\mu_{2211} & \mu_{2212} & \mu_{2221} & \mu_{2222}
\end{bmatrix}
=
\frac{1}{14}
\begin{bmatrix}
80 & 5 & 5 & 25\\
5 & 30 & 5 & 0\\
5 & 5 & 30 & 0\\
25 & 0 & 0 & 40
\end{bmatrix}.
\]

This tensor is anisotropic. Indeed, in the isotropic case, the fourth-order viscosity tensor must have the form
\[
\mu_{ijkl}^{\rm iso}
=
\alpha\,\delta_{ij}\delta_{kl}
+
\theta\,(\delta_{ik}\delta_{jl}+\delta_{il}\delta_{jk}),
\]
where $\alpha$ and $\theta$ are scalar constants. In particular, for an isotropic tensor one has
\[
\mu_{1211}^{\rm iso}=0.
\]
However, for the tensor used here,
\[
\mu_{1211}=\frac{5}{14}\neq 0.
\]
Therefore, the above tensor does not satisfy the isotropic structure and hence represents an anisotropic medium.

To generate synthetic data, the density $\rho$, the pressure $p$ and the body force $\bF$, for the numerical experiments, we use a manufactured-solution approach.
Let $\Omega=(-1,1)^2$ and $T = 0.5$. Let $\bu^0(\x)$ be the prescribed initial velocity field.
We define the space-time velocity field by
\begin{equation}
\bu(\x,t)=e^{-\lambda t}\bu^0(\x)+(1-e^{-\lambda t})\bW(\x),
\qquad \x=(x,y)\in\Omega,\ t\in(0,T),
\label{eq:manufactured_u}
\end{equation}
where $\lambda>0$ is fixed and
\begin{equation}
\mathbf{W}(x,y)=
\begin{pmatrix}
(1-x^2)(1-y^2)\\[0.2em]
(1-x^2)(1-y^2)
\end{pmatrix}.
\label{eq:manufactured_W}
\end{equation}
In our computational experiments, $\lambda = 0.2$.
Thus, $\bu(\x,0)=\bu^0(\x)$, and the added profile $\bW$ vanishes on $\partial\Omega$.

The density is initialized by a positive function $\rho_0(\x) = 1$ and is then evolved numerically from the continuity equation
\begin{equation}
\rho_t+\Div(\rho\bu)=0
\qquad \text{in } \Omega\times(0,T),
\label{eq:continuity_num}
\end{equation}
using a finite-volume update with local Lax--Friedrichs numerical fluxes in both coordinate directions.
In the experiments reported below, we set
\begin{equation}
p(\x,t)=\rho(\x,t),
\label{eq:p_equals_rho}
\end{equation}
which means $c$ in \eqref{main_eqn} is 1.

Once $\rho$ and $\bu$ are available, the body force $\bF$ is defined so that the momentum equation is satisfied numerically:
\begin{equation}
\rho\bigl(\bu_t+(\bu\cdot\nabla)\bu\bigr)
=
-\nabla p+\Div(\bmu:\nabla\bu)+\rho\bF.
\label{eq:momentum_forced_num}
\end{equation}
More precisely, $\bF$ is computed pointwise from
\begin{equation}
\bF
=
\bu_t+(\bu\cdot\nabla)\bu
+\frac{1}{\rho}\nabla p
-\frac{1}{\rho}\Div(\bmu:\nabla\bu),
\label{eq:F_definition_num}
\end{equation}
where spatial derivatives are approximated by second-order finite differences, and time derivatives are approximated by central differences in the interior and one-sided differences at the endpoints.
The anisotropic viscous term $\Div(\bmu:\nabla\bu)$ is evaluated numerically using the prescribed fourth-order tensor $\bmu$.
Therefore, the quadruple $(\rho,\bu,p,\bF)$ is consistent with the discrete forward model used to generate the boundary measurements.

After constructing the manufactured forward solution $(\rho,\bu,p,\bF)$ as described above, we extract the associated Neumann boundary data by computing the normal derivative $\partial_\nu \bu$ along the boundary $\partial\Omega$. This exact boundary flux is then used to generate synthetic measurement data by introducing multiplicative noise. More precisely, we define the noisy observations by
\[
\mathbf{f}_{\rm noise}=\mathbf{f}(1+\delta\,\mathcal R),
\qquad
\mathbf{f}=\partial_\nu\bu,
\]
where $\delta=10\%$ is the noise level and $\mathcal R$ is a uniformly distributed random field taking values in $[-1,1]$. This construction is intended to mimic measurement errors and to assess the robustness of the proposed inversion algorithm under noisy boundary observations.

\subsection{Numerical examples}

With the boundary data available, we proceed to solve the inverse problem. The main computational steps of the proposed method, based on the time-dimensional reduction framework, are summarized in Algorithm~\ref{alg}.

In Step~\ref{s1} of Algorithm~\ref{alg}, we select the truncation level $N$ and the regularization parameter $\epsilon$ through a trial-and-error calibration based on a reference test case. In the numerical experiments reported below, we choose
\[
N=20,
\qquad
\epsilon=10^{-8},
\]
and keep these values fixed across all test cases to ensure consistency and comparability across scenarios. We also use the damping parameter $\omega=0.5$ in the Picard iteration. The number of iteration $K$ in the loop Steps \ref{s4}--\ref{s7} of Algorithm \ref{alg} is 12.

In Step~\ref{s5}, the quasi-reversibility functional $J_{{\bf U}^{(k)}}^{\epsilon,N}$ is minimized to obtain an intermediate iterate $\widetilde{\bf U}^{(k+1)}$, which solves the linearized problem~\eqref{linear} in the least-squares sense. The new iterate ${\bf U}^{(k+1)}$ is then obtained through the damped Picard update introduced above. After spatial discretization, the minimization of $J_{{\bf U}^{(k)}}^{\epsilon,N}$ is carried out using standard optimization tools available in MATLAB. In our implementation, the spatial derivatives in \eqref{linear} are discretized using finite differences on a uniform $41\times 41$ Cartesian grid, and the resulting least-squares problem is solved by the built-in MATLAB function \texttt{lsqlin}. For further details on the implementation of this least-squares optimization with \texttt{lsqlin}, we refer to \cite{LeNguyen:jiip2022, Nguyen:CAMWA2020}.

The remaining steps of Algorithm~\ref{alg}, including the reconstruction of the space-time approximation $\bu^{\rm comp}$ and the recovery of the initial condition $\mathbf{U}_0^{\rm comp}$, are implemented directly from the definitions of the projected expansion.

\subsection*{Test 1.}

In the first test, the true initial velocity field
$
\bu^0_{\rm true}
=
\begin{bmatrix}
u^0_{{\rm true},1} & u^0_{{\rm true},2}
\end{bmatrix}^{\top}
$
is chosen so that its two components are supported in distinct elliptical regions. Specifically, we set
\[
u^0_{{\rm true},1}(x,y)=
\begin{cases}
1, & \text{if } 8(x-0.4)^2+y^2<0.4^2,\\
0, & \text{otherwise},
\end{cases}
\qquad
u^0_{{\rm true},2}(x,y)=
\begin{cases}
1, & \text{if } x^2+8(y-0.3)^2<0.4^2,\\
0, & \text{otherwise}.
\end{cases}
\]
The function $u^0_{{\rm true},1}$ is supported in a vertically elongated ellipse centered near $(0.4,0)$, while $u^0_{{\rm true},2}$ is supported in a horizontally elongated ellipse centered near $(0,0.3)$. This choice allows us to examine the ability of the proposed method to recover both the location and the geometry of localized anisotropic structures in the initial data. The boundary observations are generated synthetically from the manufactured data $(\rho,\bu,p,\bF)$ and then used as input for the inverse reconstruction algorithm.

\begin{figure}[h!]
\centering
	\subfloat[\label{test1a} True first velocity component $u^0_{\mathrm{true},1}$]{\includegraphics[width=.25\textwidth]{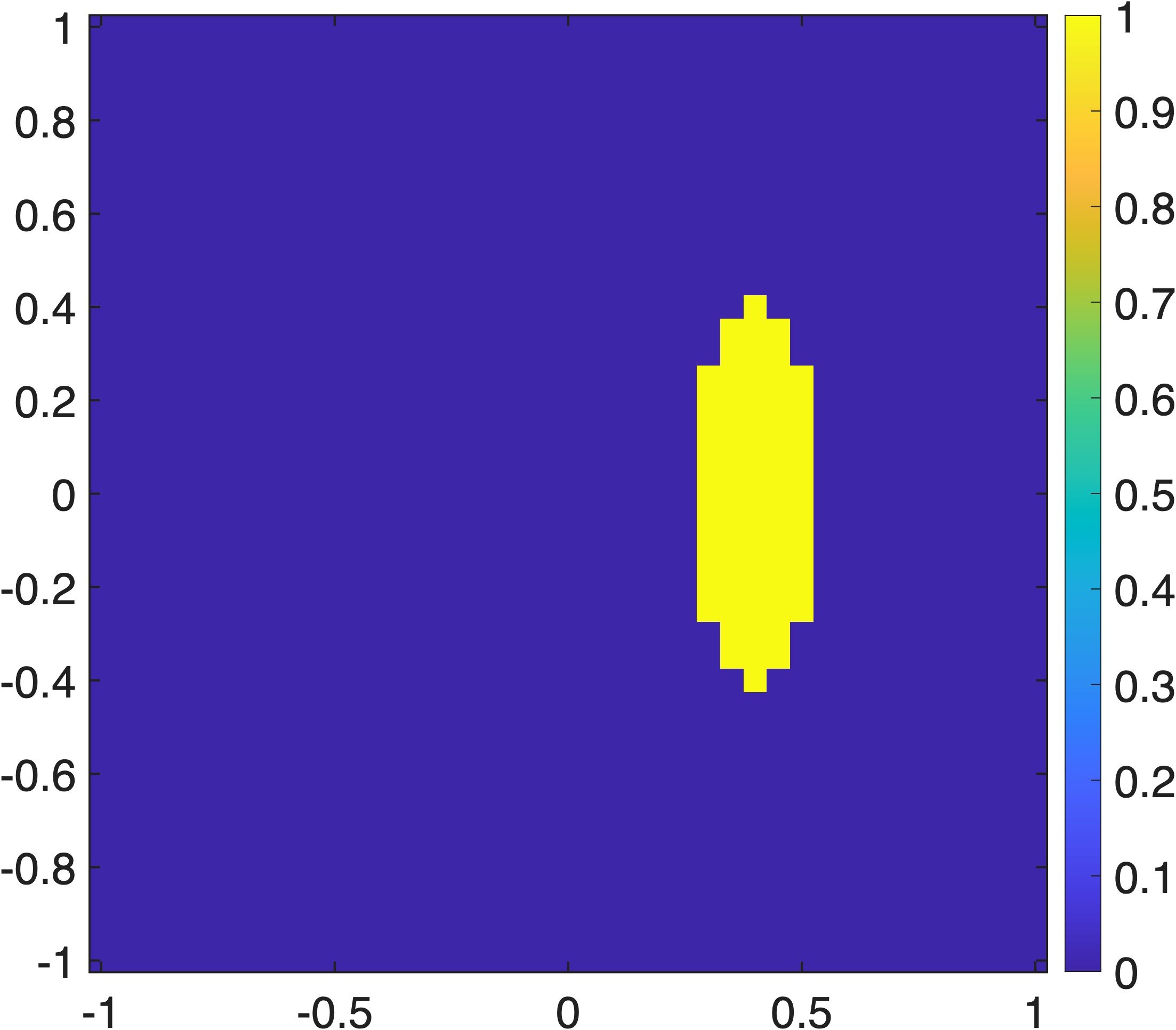}}
	\quad
	\subfloat[\label{test1b} Reconstructed first velocity component $u^0_{\mathrm{comp},1}$]{\includegraphics[width=.25\textwidth]{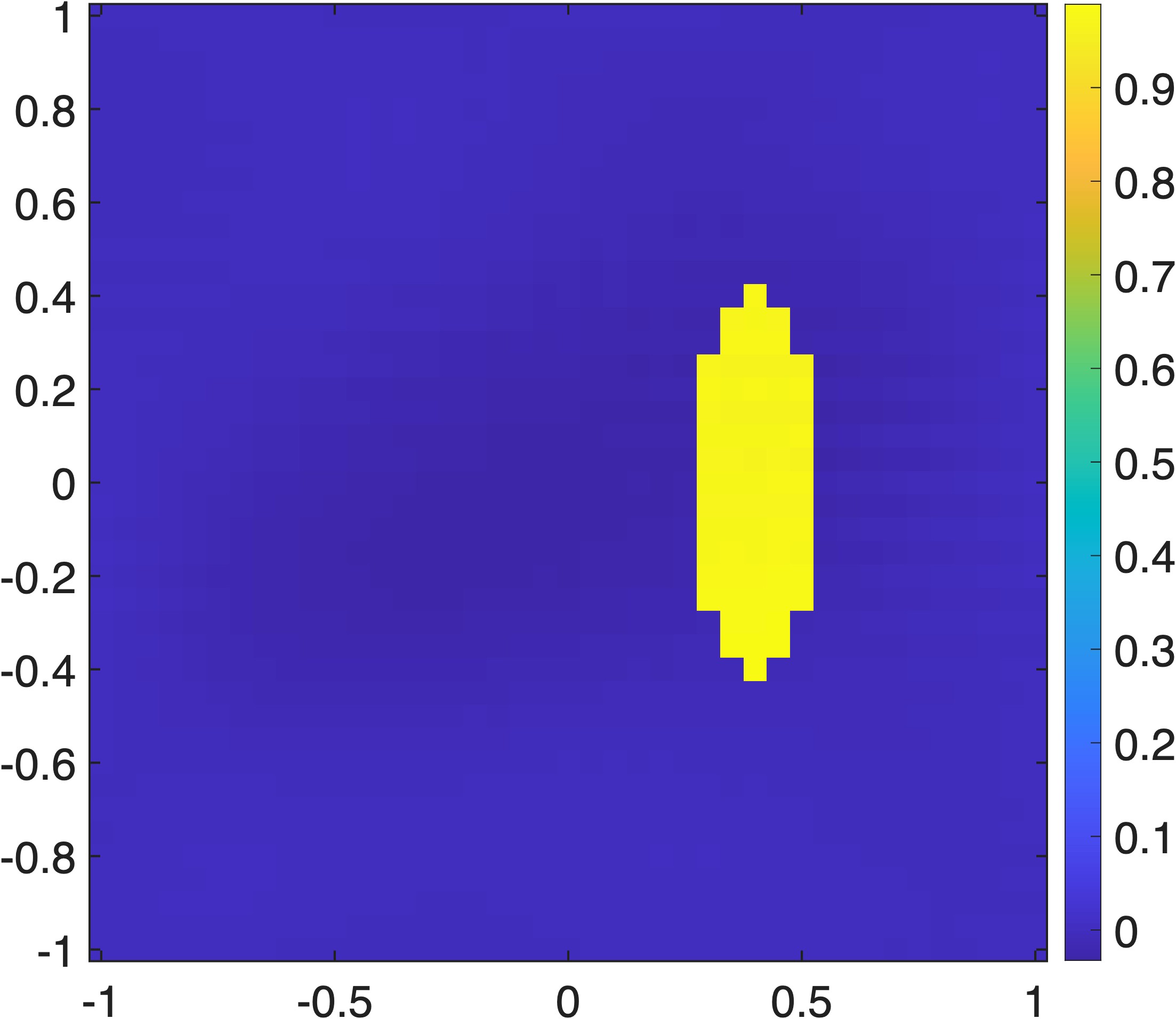}}
	\quad
	\subfloat[\label{test1c} Relative error of the first component $\frac{|u^0_{{\rm true},1} - u^0_{{\rm comp},1}|}{\|u^0_{{\rm true},1}\|_{L^\infty(\Omega)}}$
]{\includegraphics[width=.25\textwidth]{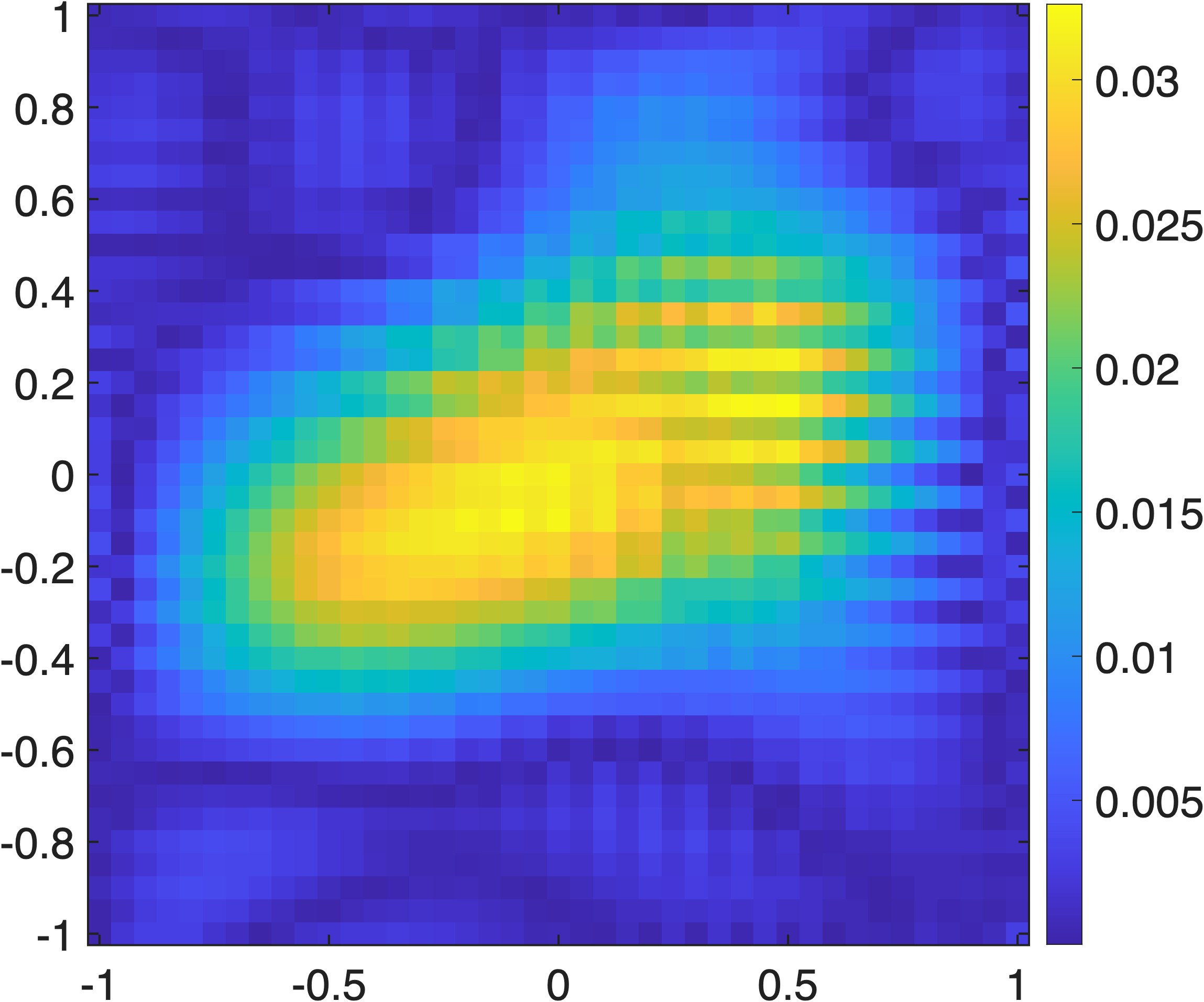}}
	
	\subfloat[\label{test1d}  True second velocity component $u^0_{\mathrm{true},2}$]{\includegraphics[width=.25\textwidth]{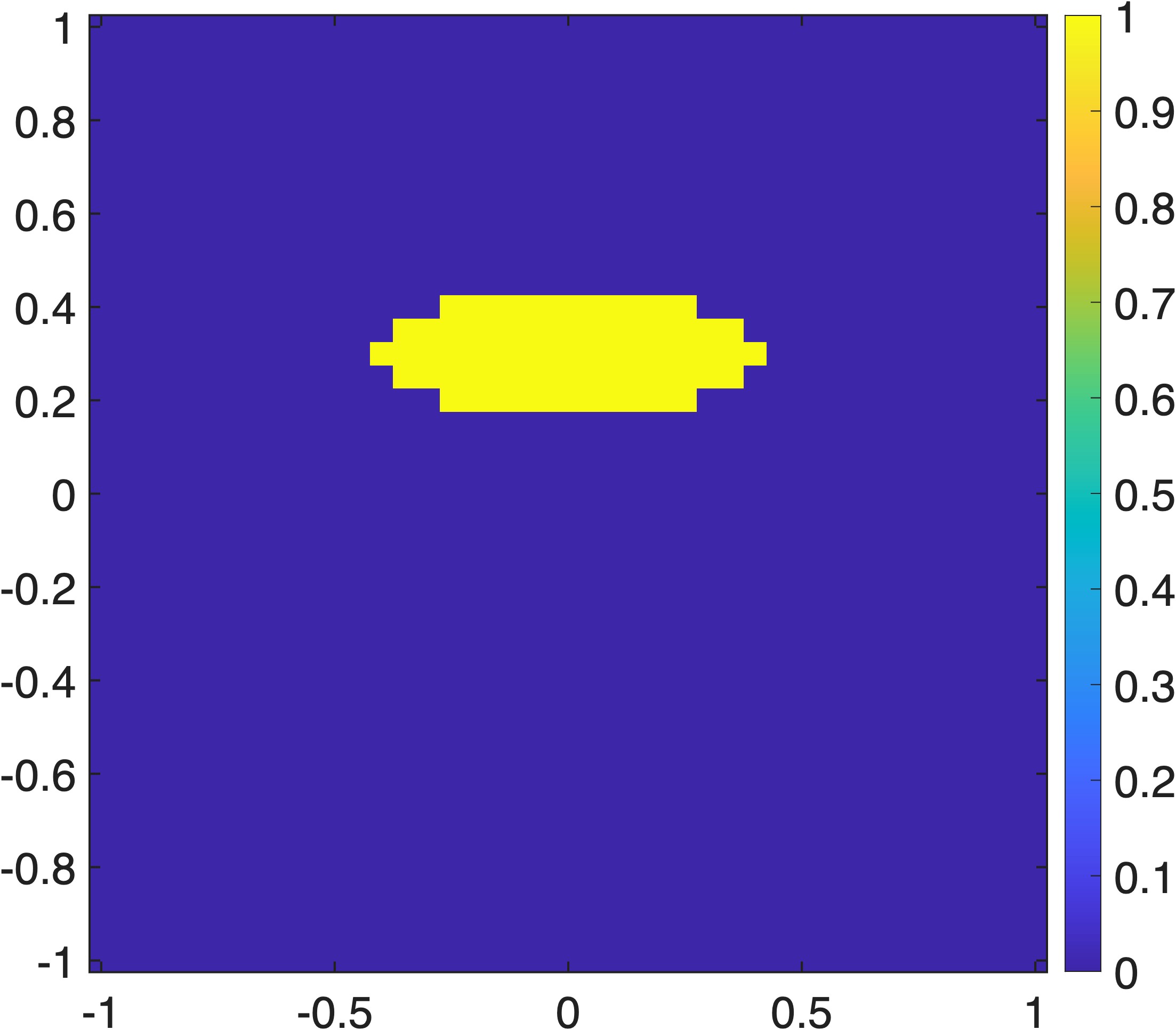}}
	\quad
	\subfloat[\label{test1e} Reconstructed second velocity component $u^0_{\mathrm{comp},2}$]{\includegraphics[width=.25\textwidth]{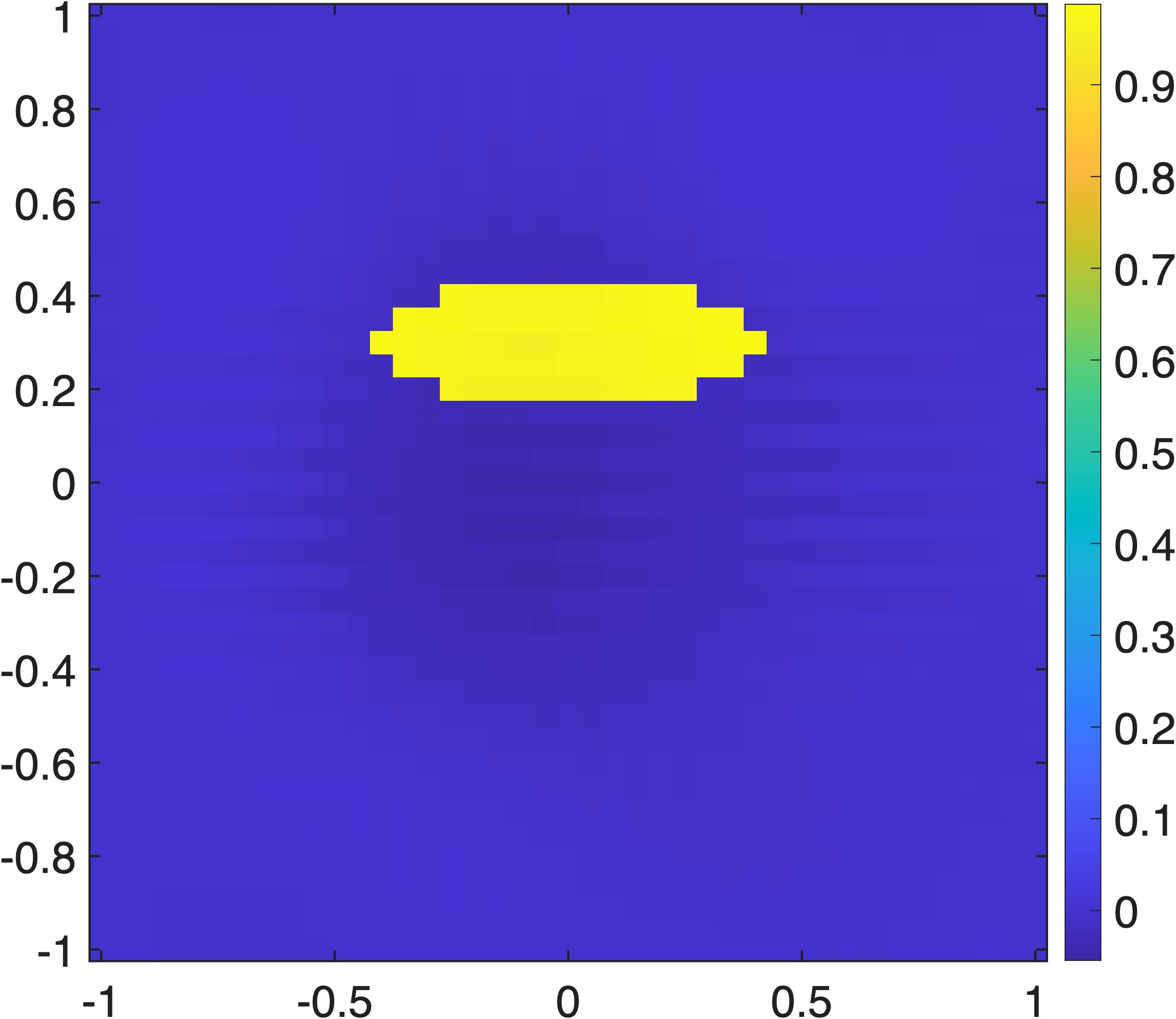}}
	\quad
	\subfloat[\label{test1f} Relative error of the second component $\frac{|u^0_{{\rm true},2} - u^0_{{\rm comp},1}|}{\|u^0_{{\rm true},2}\|_{L^\infty(\Omega)}}$]{\includegraphics[width=.25\textwidth]{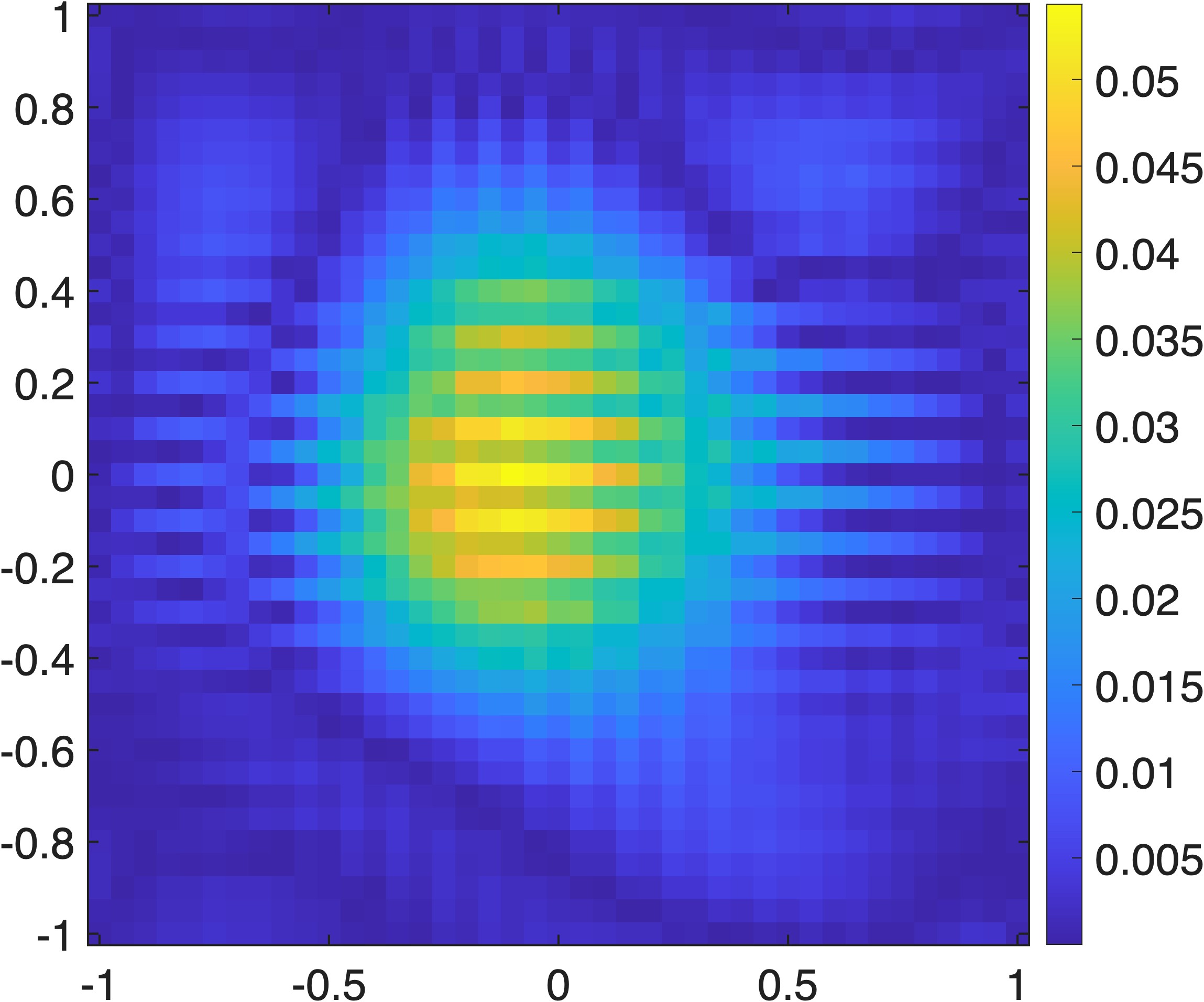}}
	
	\subfloat[\label{test1g}  True initial velocity field $\mathbf{u}^0_{\mathrm{true}}$]{\includegraphics[width=.25\textwidth]{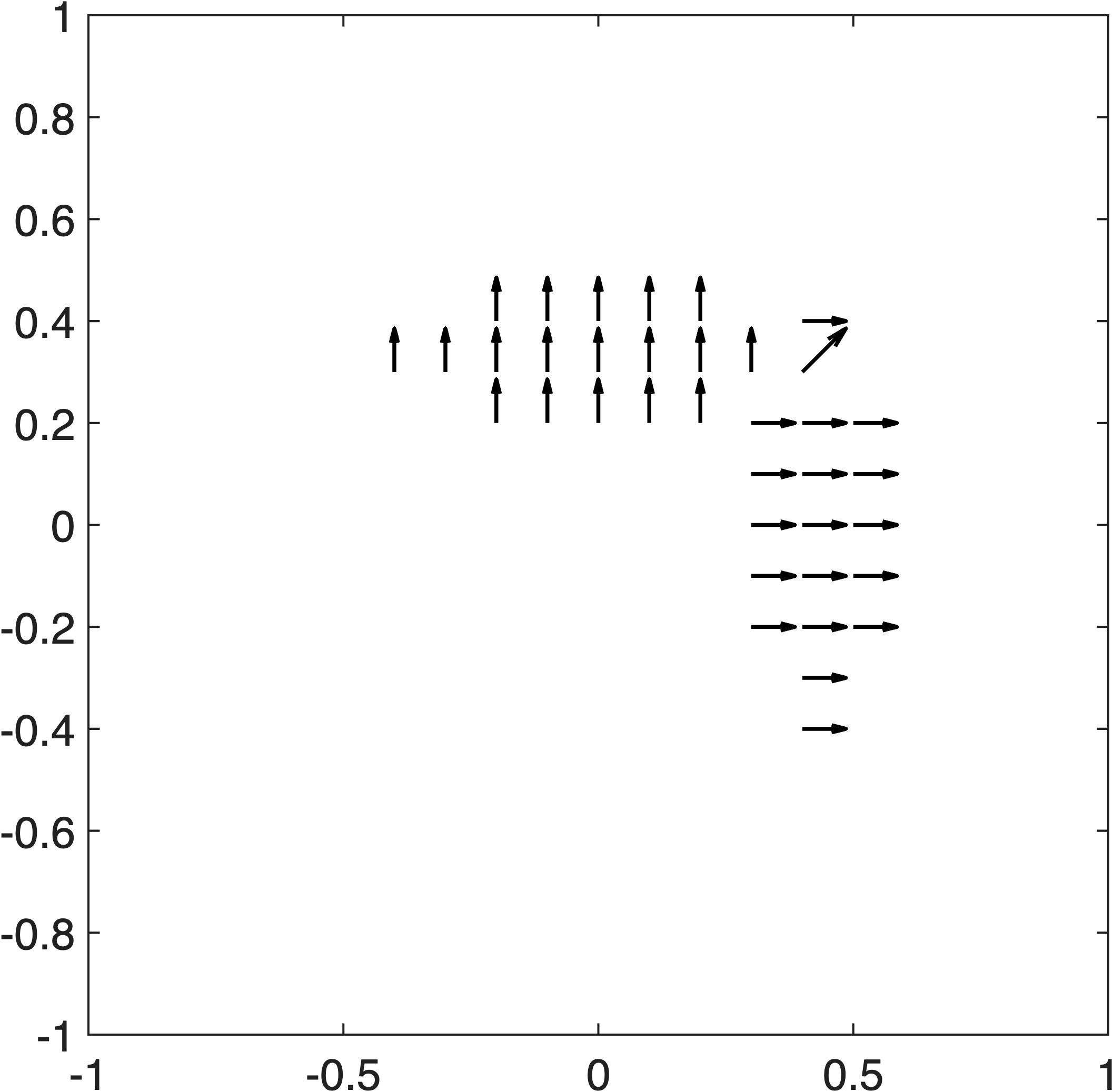}}
	\quad
	\subfloat[\label{test1h} Reconstructed initial velocity field $\mathbf{u}^{0}_{\mathrm{comp}}$]{\includegraphics[width=.25\textwidth]{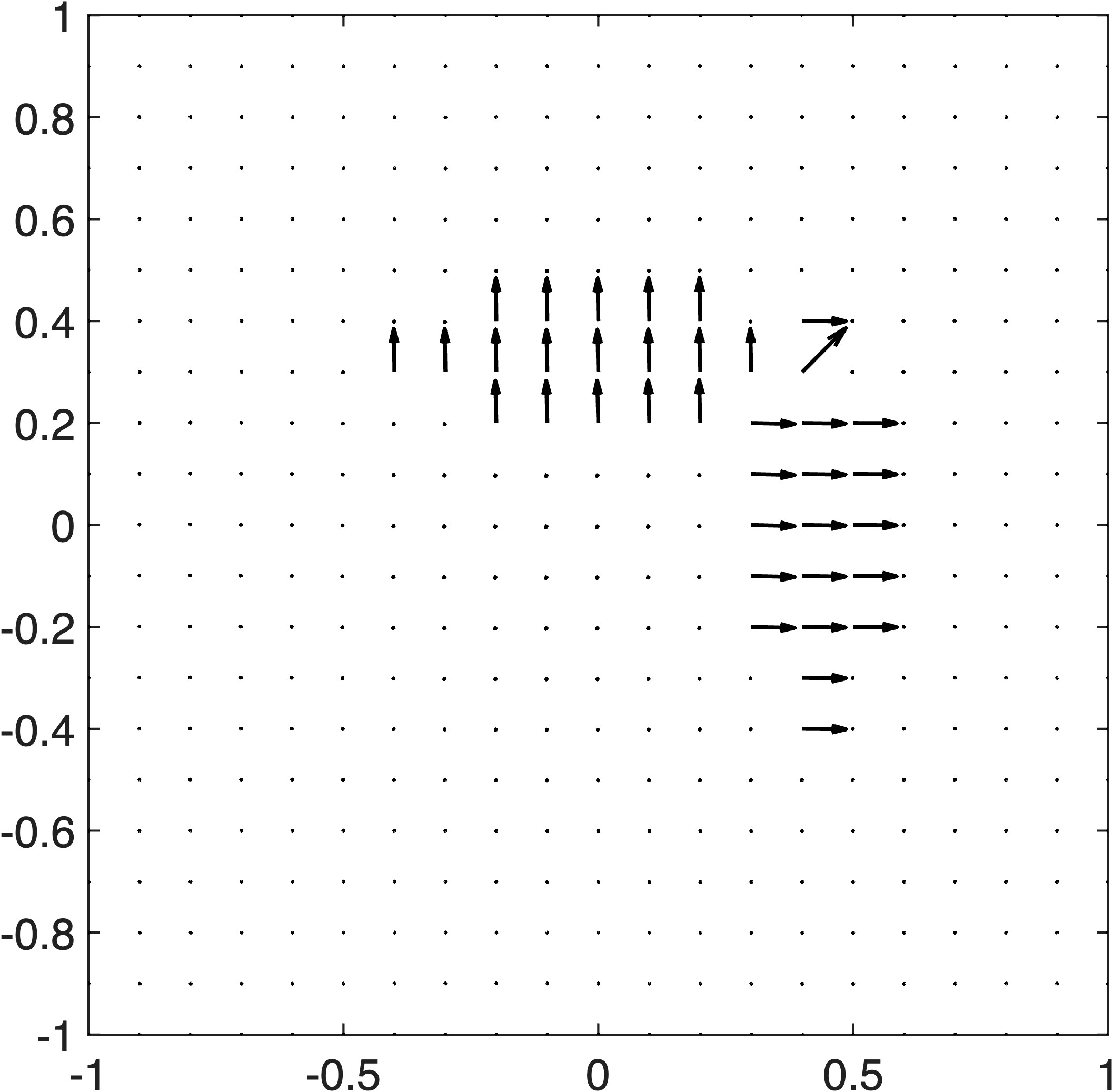}}
	\quad
	\subfloat[\label{test1i} Residual of \eqref{time_red} for $\bU^{(k)}$ generated by Algorithm \ref{alg}]{\includegraphics[width=.25\textwidth]{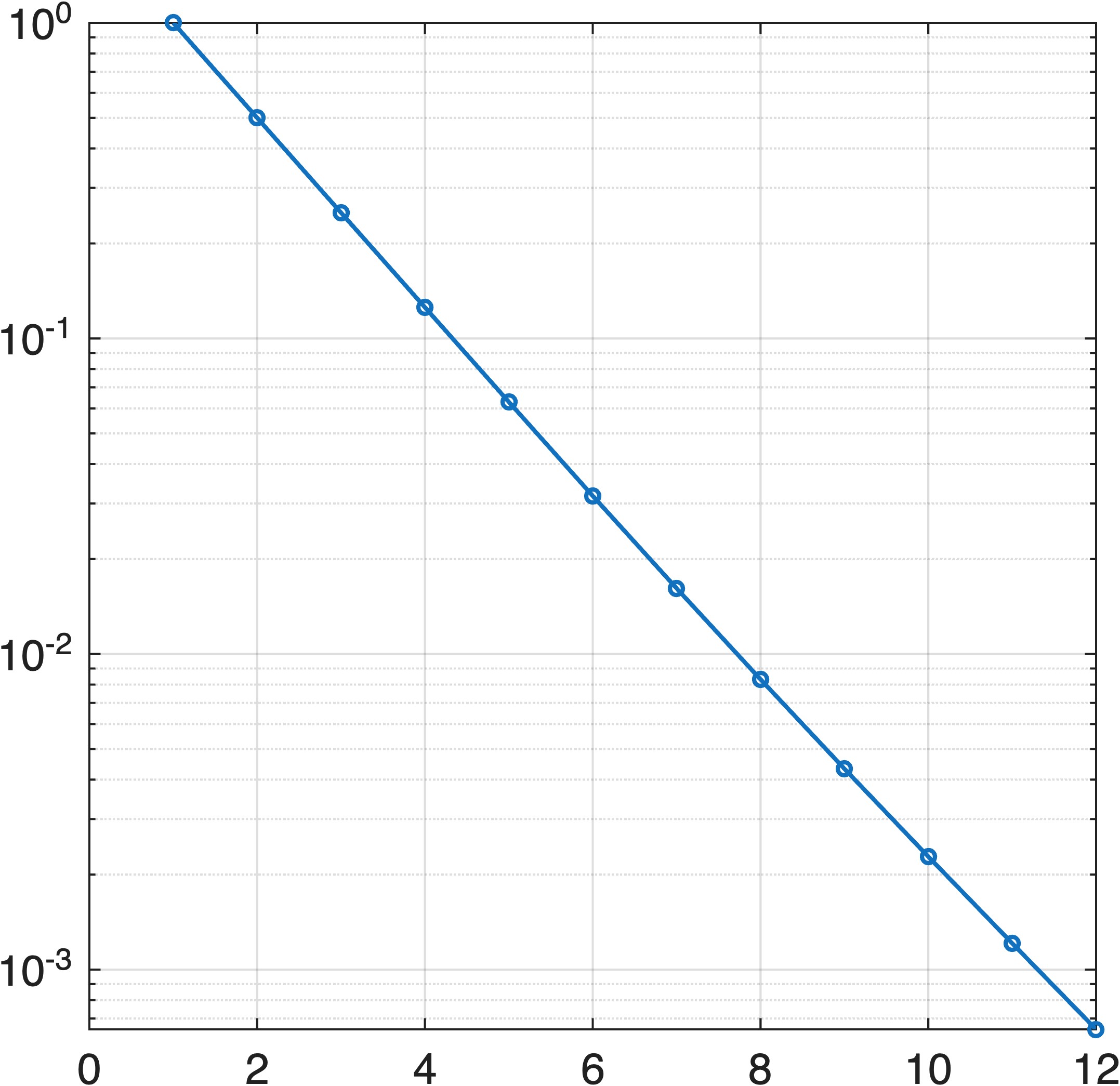}}

	\caption{\label{test1} Test 1: true and reconstructed initial velocity fields. Figures (a)--(b) compare the true and reconstructed first velocity component, while Figure (c) displays the corresponding relative error. Figures (d)--(e) compare the true and reconstructed second velocity component, while Figure (f) shows the corresponding relative error. Figures (g)--(h) present the true and reconstructed initial velocity fields, respectively. Figure (i) plots the residual of \eqref{time_red} for $\bU^{(k)}$ generated by Algorithm~\ref{alg}, showing that $\bU^{(k)}$ approximately solves the time-reduction model \eqref{time_red}.}
\end{figure}

The results presented in Figure~\ref{test1} demonstrate the strong performance of the proposed method for the inverse initial data problem, even in the presence of $10\%$ noise in the boundary measurements. Figures~\ref{test1a}--\ref{test1b} and Figures~\ref{test1d}--\ref{test1e} show that the reconstructed first and second velocity components recover not only the correct location and support of the true solution, but also its essential geometric structure. Although the reconstructions exhibit slight diffusion and mild background artifacts, the main features of both components are captured accurately. This is further confirmed by the relative error plots in Figures~\ref{test1c} and \ref{test1f}, where the displayed maxima are approximately $0.035$ and $0.055$, respectively. In addition, the comparison between Figures~\ref{test1g} and \ref{test1h} shows that the reconstructed initial velocity field preserves the dominant directional pattern of the true field remarkably well. Finally, Figure~\ref{test1i} reveals a clear and steady decay of the residual of \eqref{time_red} for $\bU^{(k)}$, providing strong numerical evidence for the stability and effectiveness of the iterative reconstruction procedure. Altogether, these results highlight the robustness of the proposed approach and underscore its strong potential for solving noisy inverse problems in more realistic and challenging settings.

\subsection*{Test 2}
For the second test case, the true initial velocity field comprises two components, each with localized features in the form of diagonal ellipses. Specifically, the first component, $u^0_{\rm true,1}$, is defined as:
\[
u^0_{\rm true,1}(x, y) = 
\begin{cases}
1, & \text{if } 1.5(x + y)^2 + y^2 < 0.5^2, \\
0, & \text{otherwise}.
\end{cases}
\]
The second component, $u^0_{\rm true,2}$, takes the form:
\[
u^0_{\rm true,2}(x, y) = 
\begin{cases}
1, & \text{if } x^2 + 1.5(x - y)^2 < 0.5^2, \\
0, & \text{otherwise}.
\end{cases}
\]
Both components feature diagonal ellipses with distinct orientations, which serve as a test for the method's ability to recover such geometrically oriented structures.

\begin{figure}[h!]
\centering
	\subfloat[\label{test2a}True first velocity component $u^0_{\mathrm{true},1}$]{\includegraphics[width=.25\textwidth]{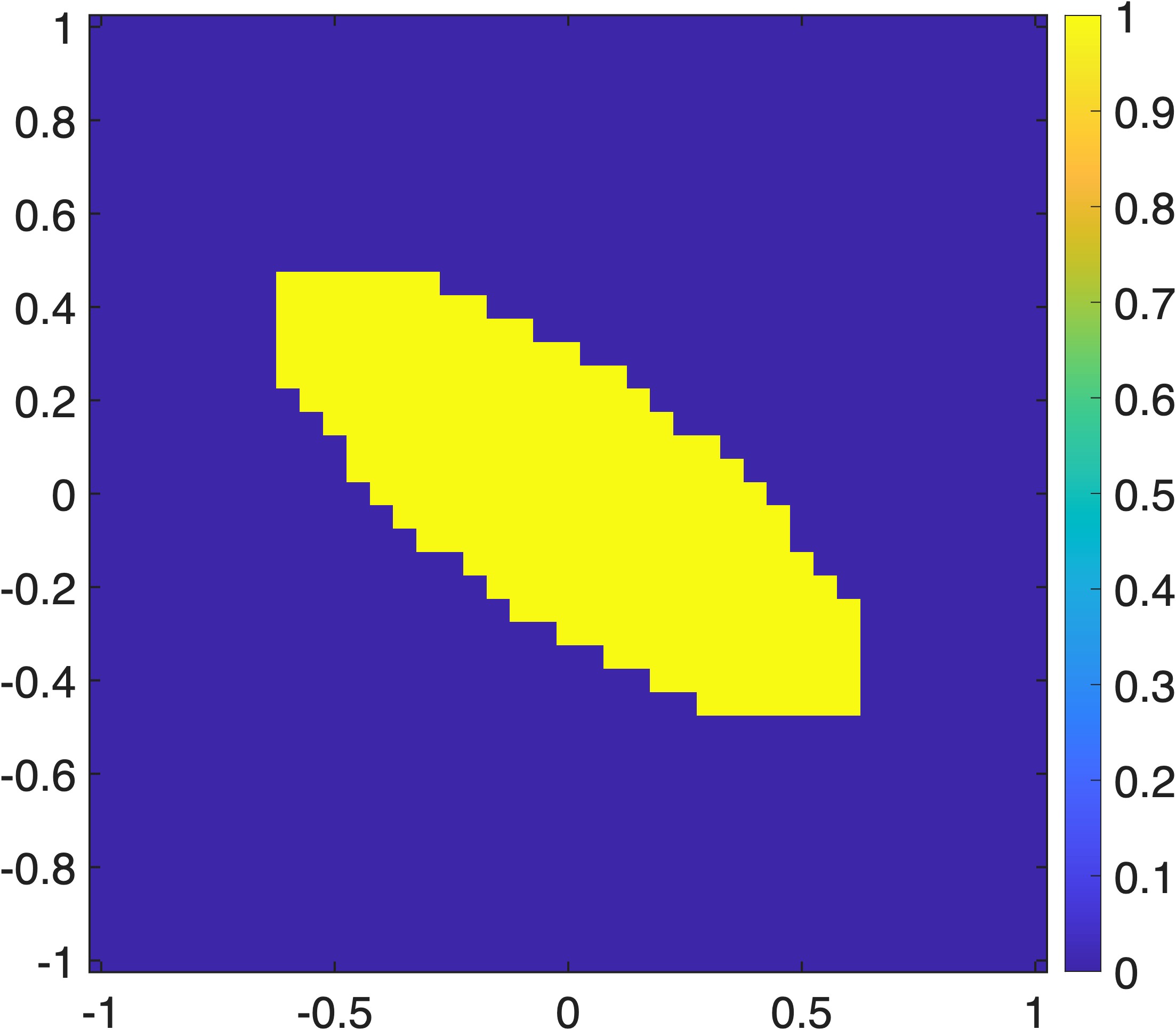}}
	\quad
	\subfloat[\label{test2b}Reconstructed first velocity component $u^0_{\mathrm{comp},1}$]{\includegraphics[width=.25\textwidth]{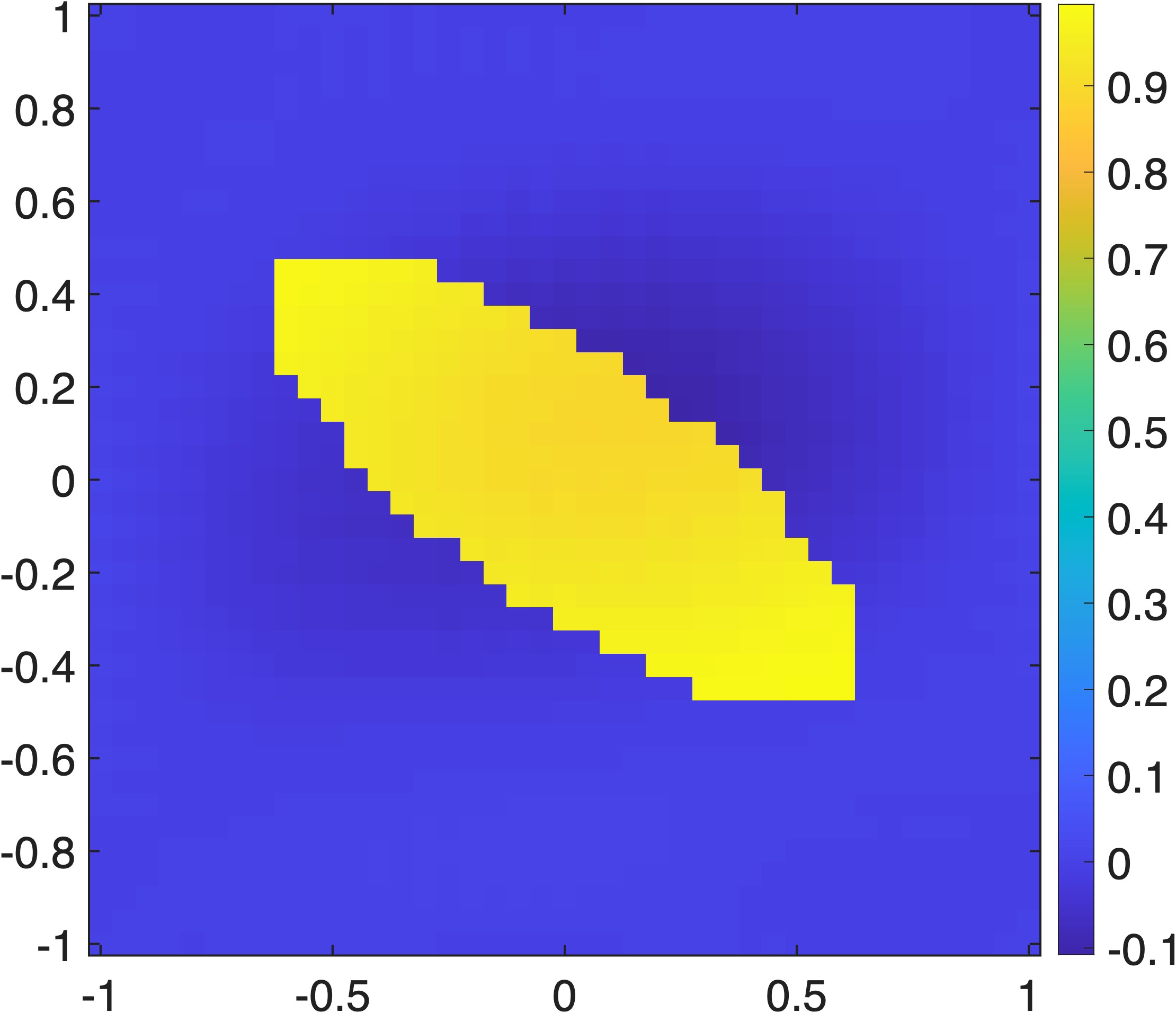}}
	\quad
	\subfloat[\label{test2c}Relative error of the first component $\frac{|u^0_{{\rm true},1} - u^0_{{\rm comp},1}|}{\|u^0_{{\rm true},1}\|_{L^\infty(\Omega)}}$
]{\includegraphics[width=.25\textwidth]{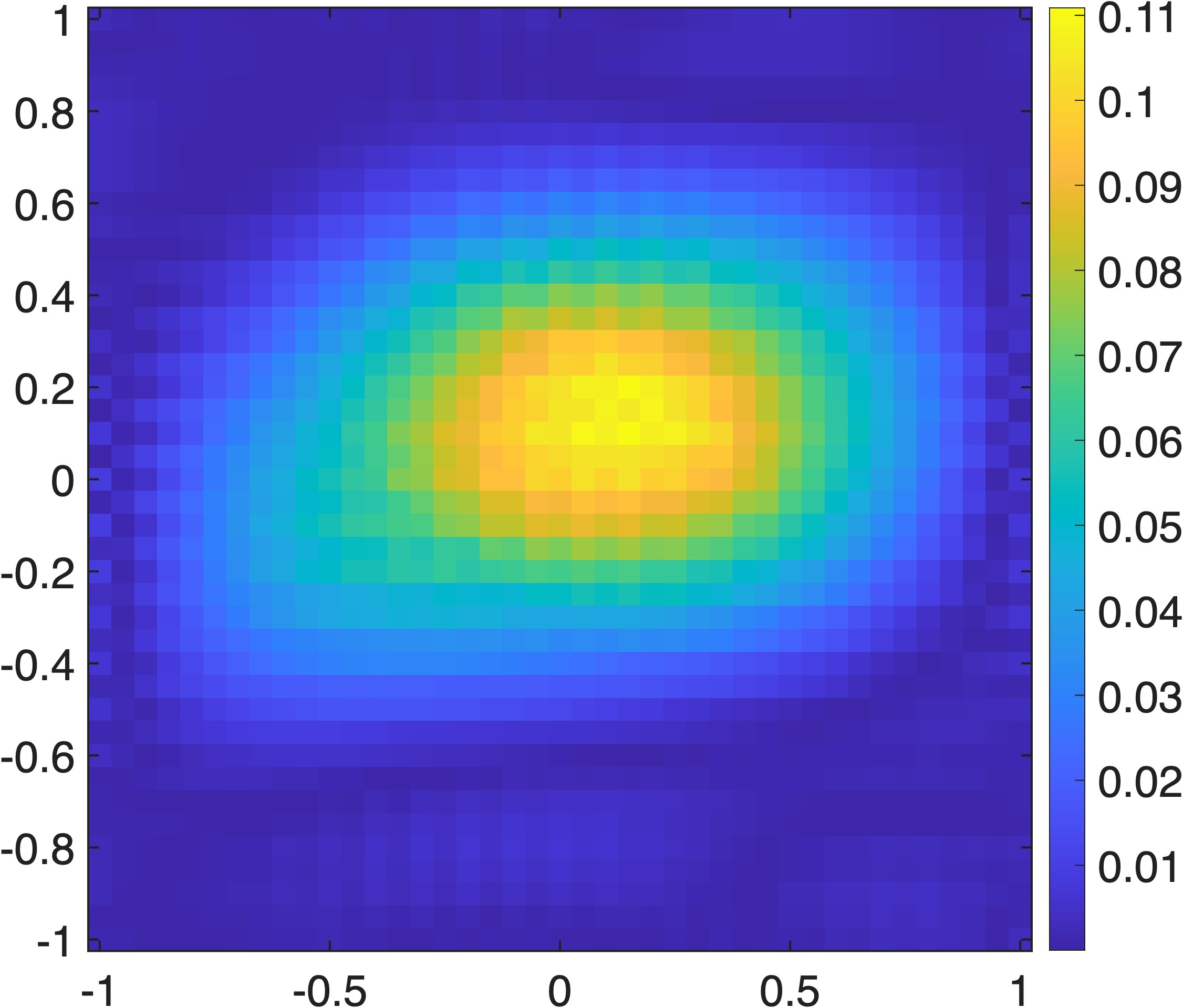}}
	
	\subfloat[\label{test2d} True second velocity component $u^0_{\mathrm{true},2}$]{\includegraphics[width=.25\textwidth]{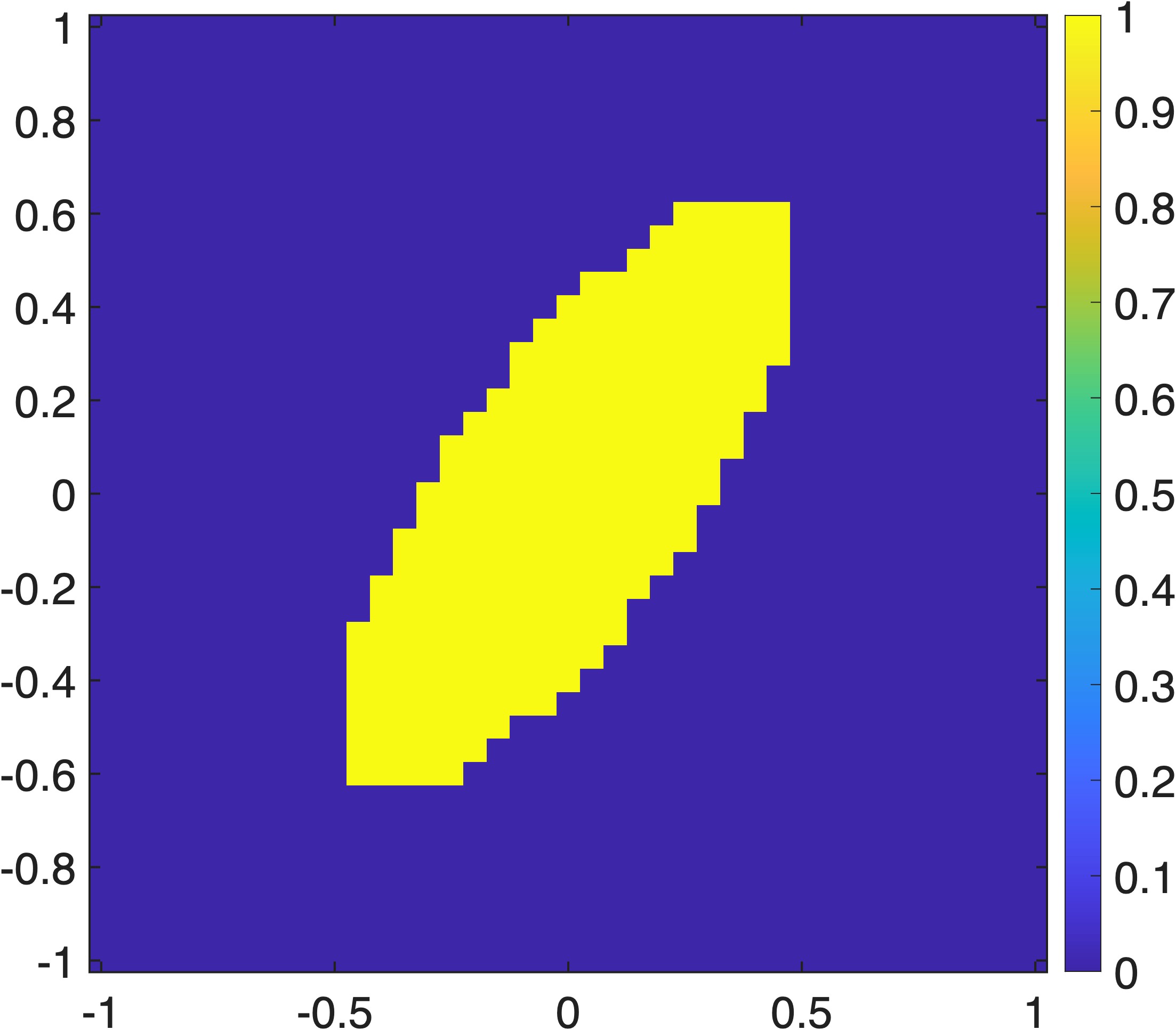}}
	\quad
	\subfloat[\label{test2e}Reconstructed second velocity component $u^0_{\mathrm{comp},2}$]{\includegraphics[width=.25\textwidth]{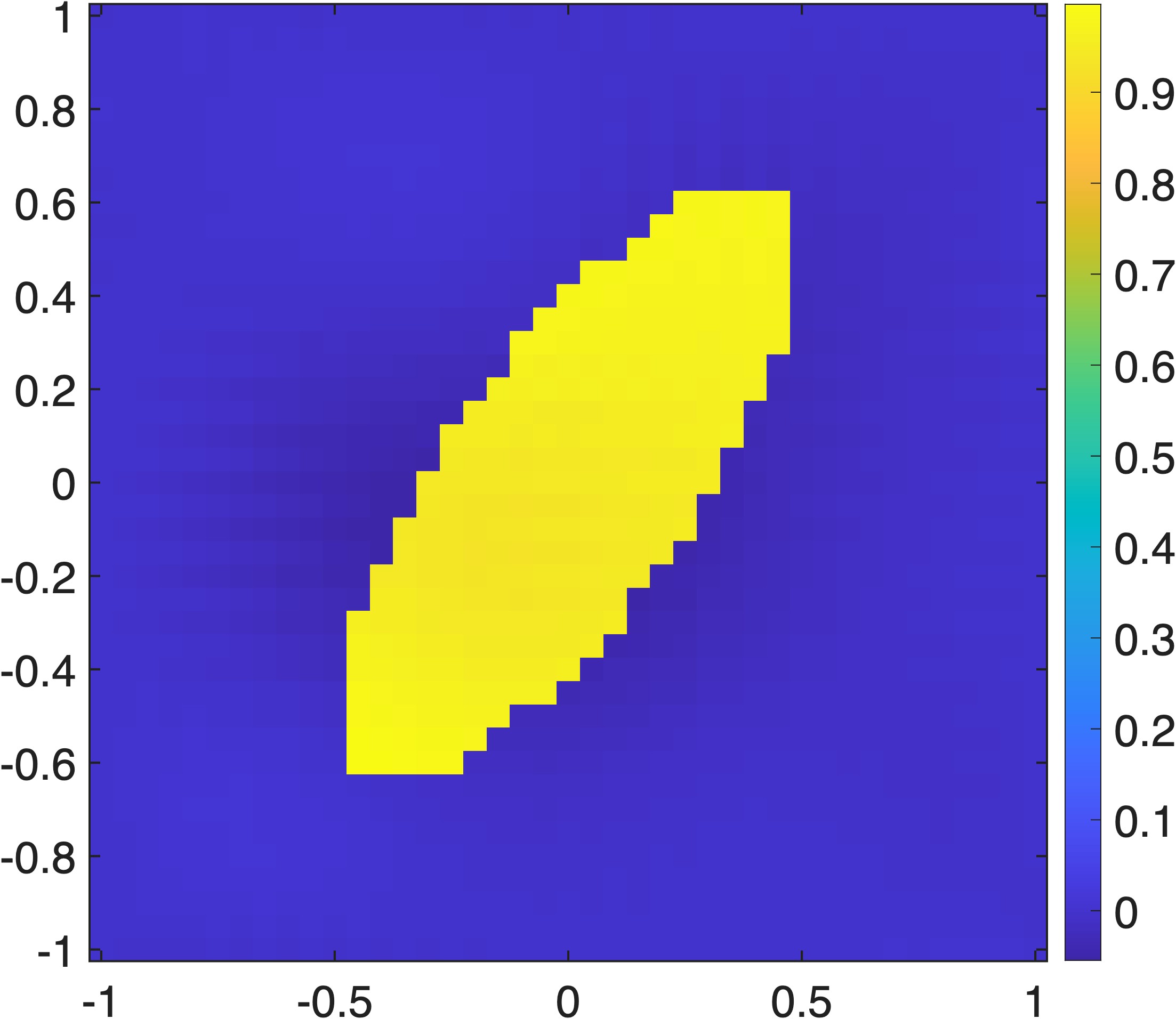}}
	\quad
	\subfloat[\label{test2f}Relative error of the second component $\frac{|u^0_{{\rm true},2} - u^0_{{\rm comp},2}|}{\|u^0_{{\rm true},2}\|_{L^\infty(\Omega)}}$]{\includegraphics[width=.25\textwidth]{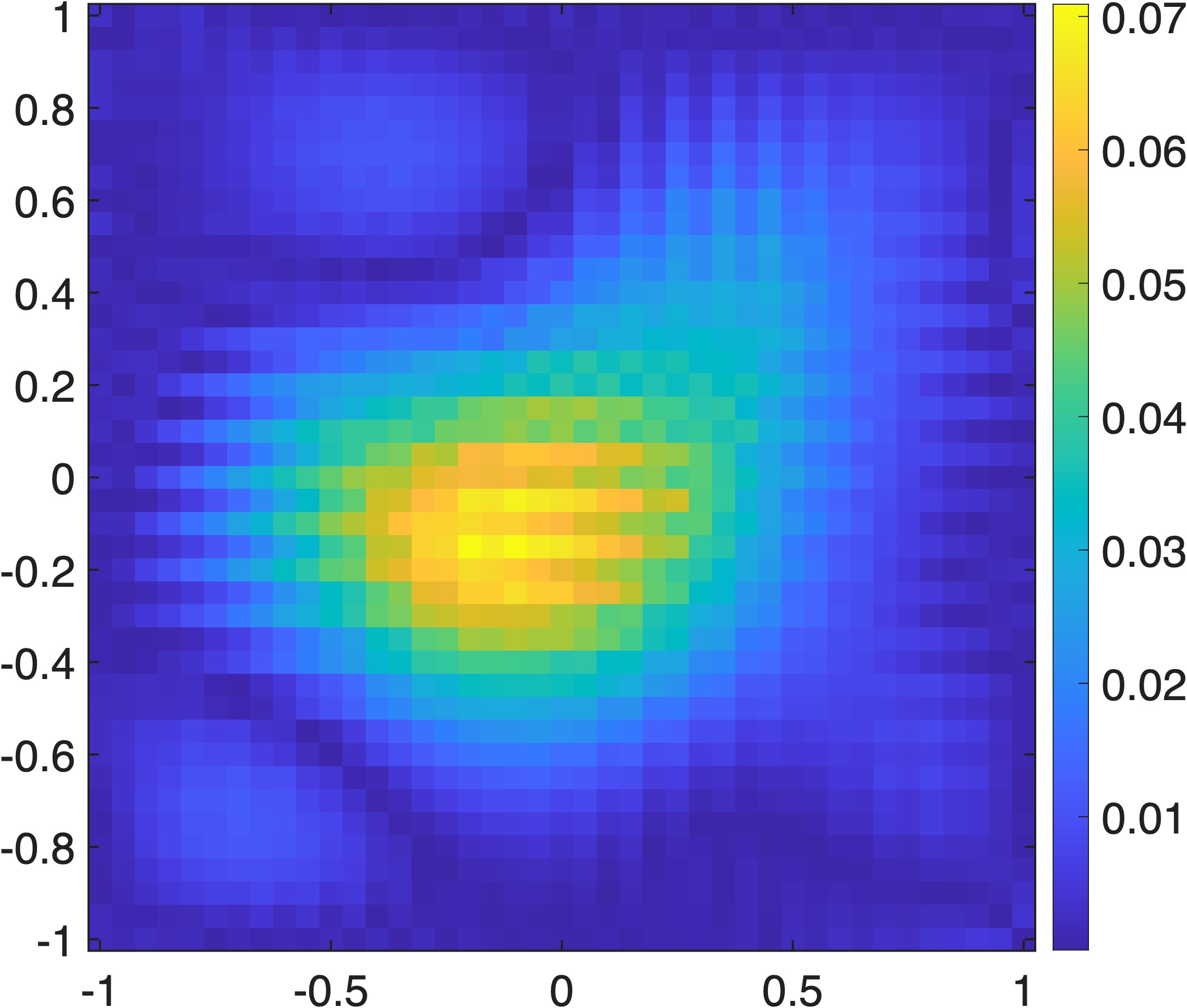}}
	
	\subfloat[\label{test2g} True initial velocity field $\mathbf{u}^0_{\mathrm{true}}$]{\includegraphics[width=.25\textwidth]{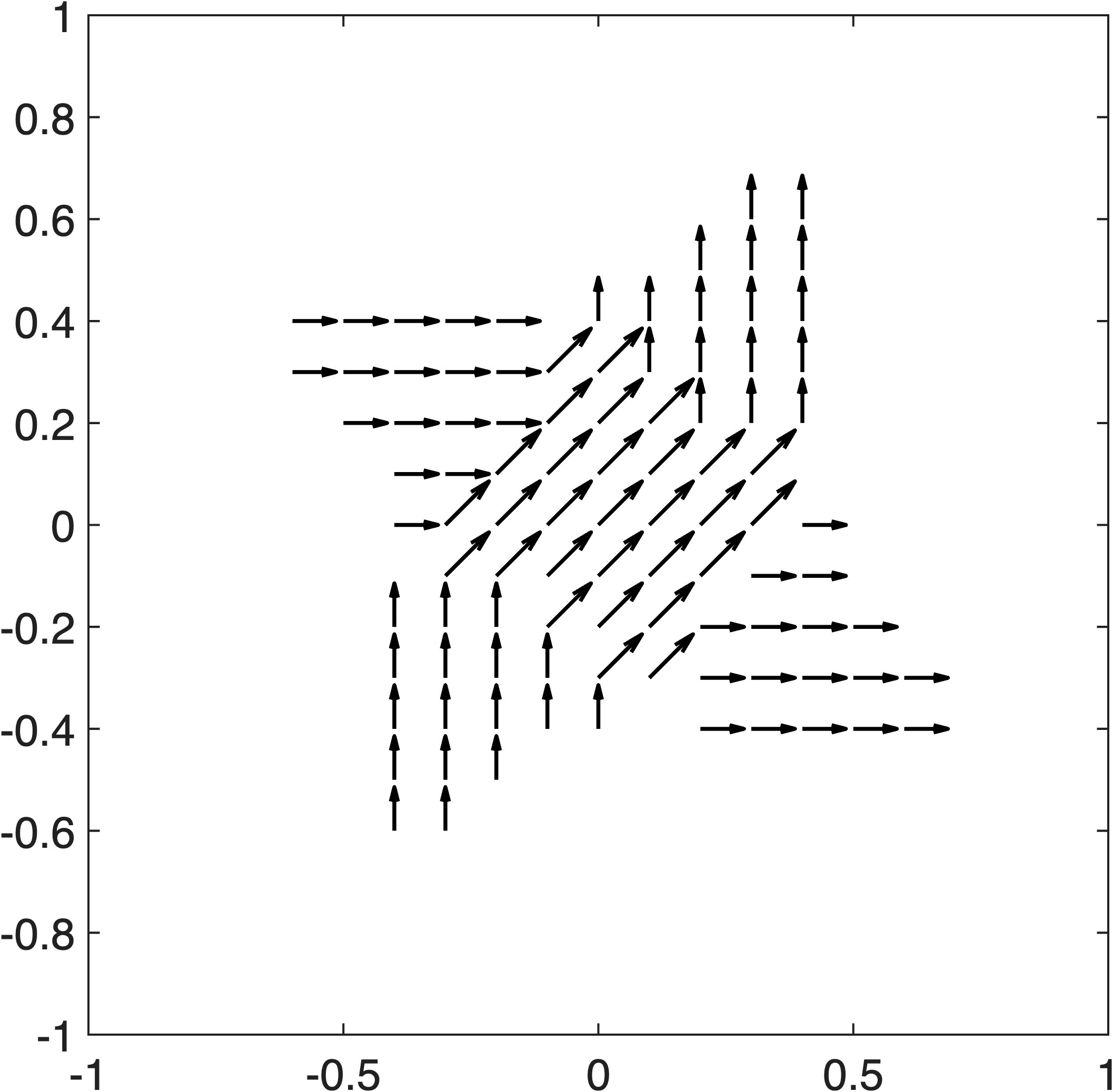}}
	\quad
	\subfloat[\label{test2h}Reconstructed initial velocity field $\mathbf{u}^{0}_{\mathrm{comp}}$]{\includegraphics[width=.25\textwidth]{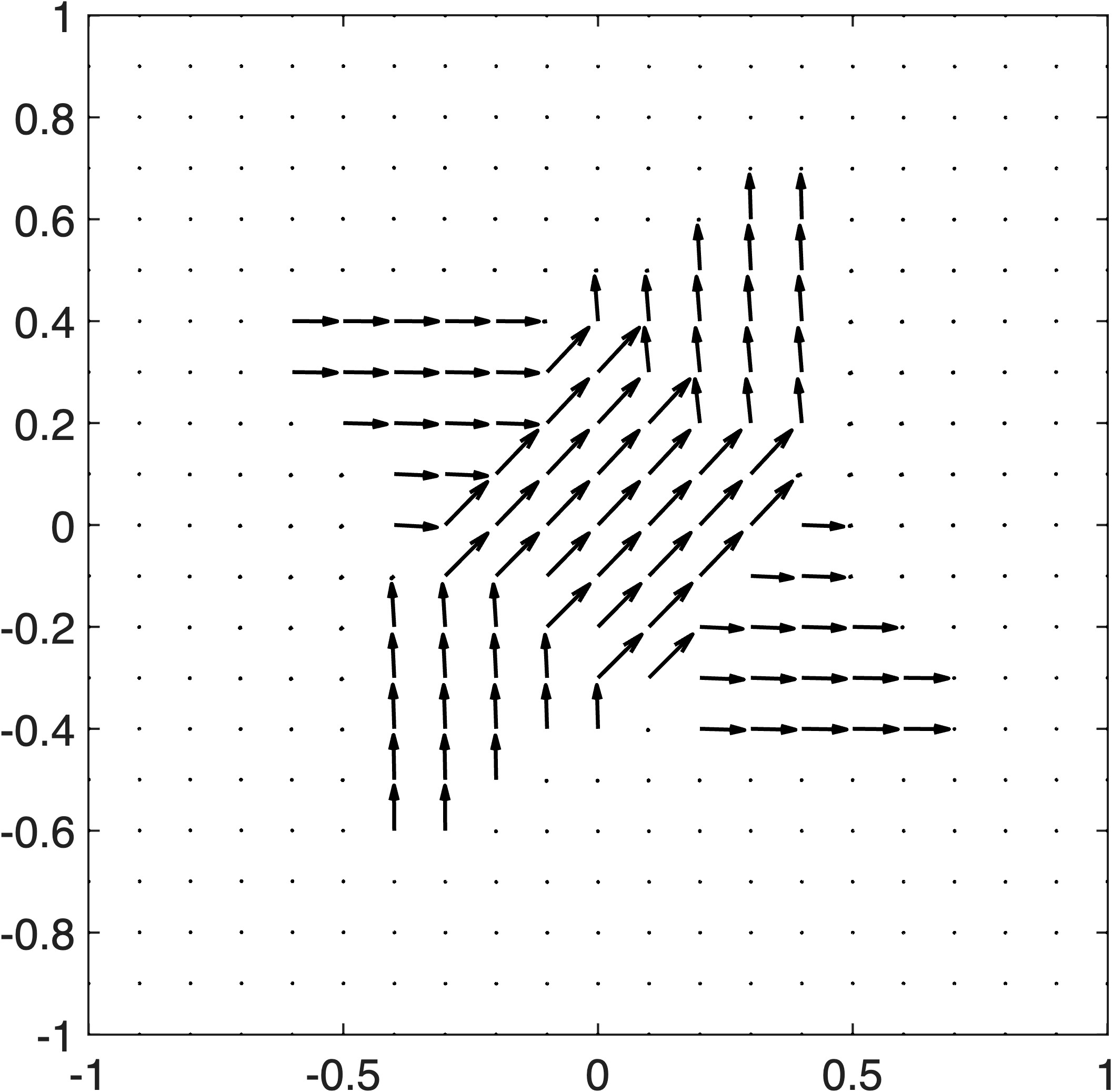}}
	\quad
	\subfloat[\label{test2i}Residual of \eqref{time_red} for $\bU^{(k)}$ generated by Algorithm \ref{alg}]{\includegraphics[width=.25\textwidth]{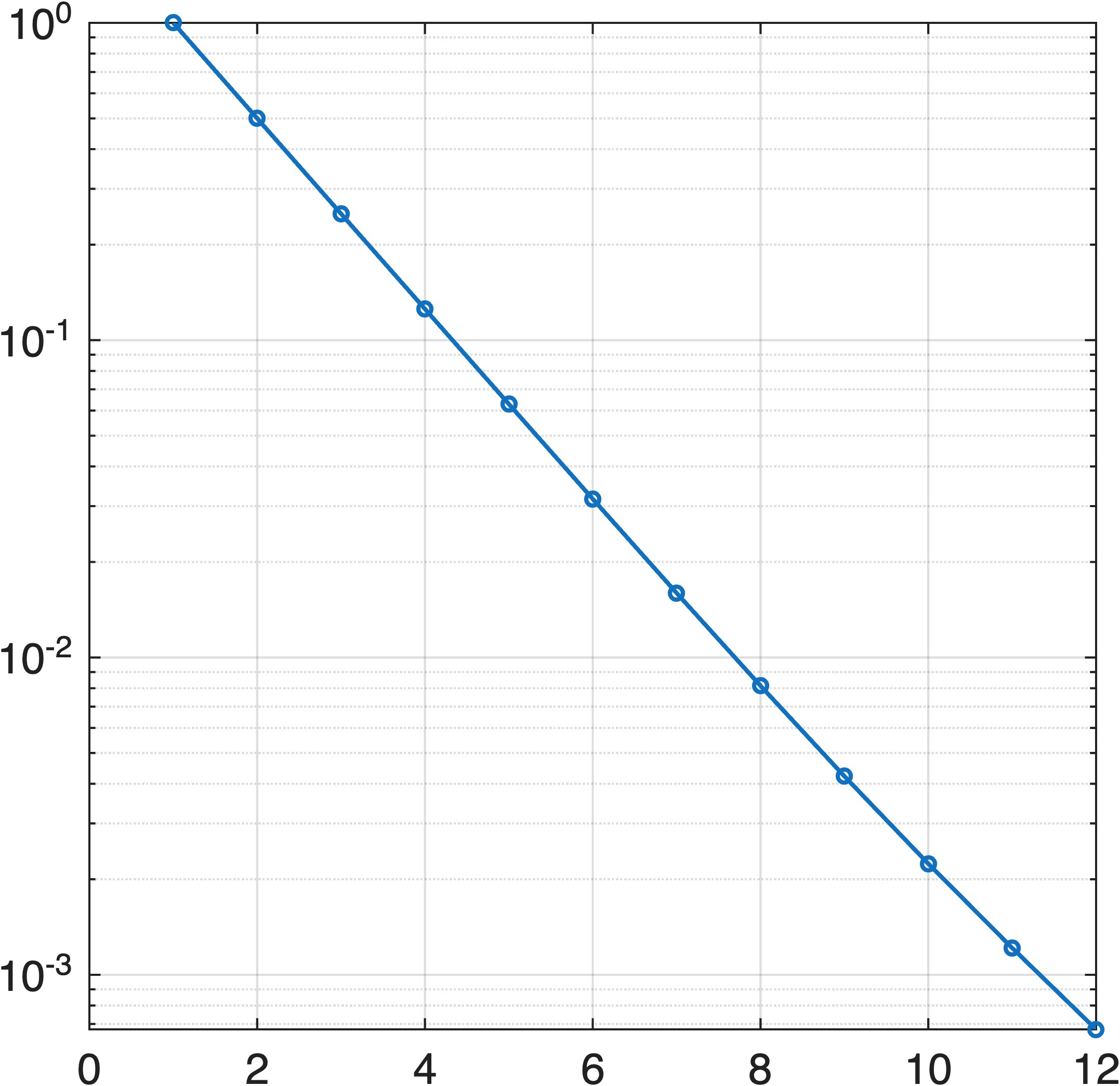}}

	\caption{\label{test2} Test 2: true and reconstructed initial velocity fields. Figures (a)--(c) show the true first component, its reconstruction, and the corresponding relative error. Figures (d)--(f) show the true second component, its reconstruction, and the corresponding relative error. Figures (g)--(h) display the true and reconstructed initial velocity fields. Figure (i) shows the residual of \eqref{time_red} for $\bU^{(k)}$ generated by Algorithm~\ref{alg}, which decreases steadily throughout the iteration. }
\end{figure}

The results shown in Figure~\ref{test2} further demonstrate the effectiveness of the proposed method in reconstructing initial velocity fields with non-axis-aligned geometric structures. Figure~\ref{test2a} and Figure~\ref{test2b}, as well as Figure~\ref{test2d} and Figure~\ref{test2e}, indicate that both diagonal components are recovered with the correct orientation, location, and overall support. Although the reconstructed profiles are slightly smoothed and accompanied by mild background artifacts, the dominant diagonal features remain clearly visible. This is supported by the relative error distributions in Figure~\ref{test2c} and Figure~\ref{test2f}, whose displayed maxima are approximately $0.11$ and $0.07$, respectively. Moreover, the comparison between Figure~\ref{test2g} and Figure~\ref{test2h} shows that the reconstructed initial velocity field preserves the main directional structure of the true field. Finally, Figure~\ref{test2i} shows a steady decay of the residual of \eqref{time_red} for $\bU^{(k)}$, which again provides numerical evidence for the stability and reliability of the damped Picard iteration. Overall, this test confirms that the proposed method is capable of accurately recovering localized anisotropic structures with different orientations, even under noisy boundary observations.

\subsection*{Test 3} 
In the third test, the true initial velocity field is chosen so that both components contain a multi-layer structure with an outer ring and an inner core of different magnitude. The first component, $u^0_{\rm true,1}(x,y)$, is defined by
\[
u^0_{\rm true,1}(x,y)=
\begin{cases}
1, & \text{if } 0.6^2<x^2+y^2<0.7^2,\\
2, & \text{if } x^2+y^2<0.2^2,\\
0, & \text{otherwise},
\end{cases}
\]
so that it consists of a circular outer annulus together with a central disk of larger positive magnitude. The second component, $u^0_{\rm true,2}(x,y)$, is defined by
\[
u^0_{\rm true,2}(x,y)=
\begin{cases}
1, & \text{if } 0.6^2<\max(x^2,y^2)<0.7^2,\\
-2, & \text{if } \max(x^2,y^2)<0.2^2,\\
0, & \text{otherwise},
\end{cases}
\]
which represents a square outer ring together with a central square core of negative magnitude. This test is designed to assess whether the proposed method can recover composite structures with multiple active regions, sharp interfaces, and sign changes across different components of the initial velocity field.
\begin{figure}[h!]
\centering
	\subfloat[\label{fig3a}True first velocity component $u^0_{\mathrm{true},1}$]{\includegraphics[width=.25\textwidth]{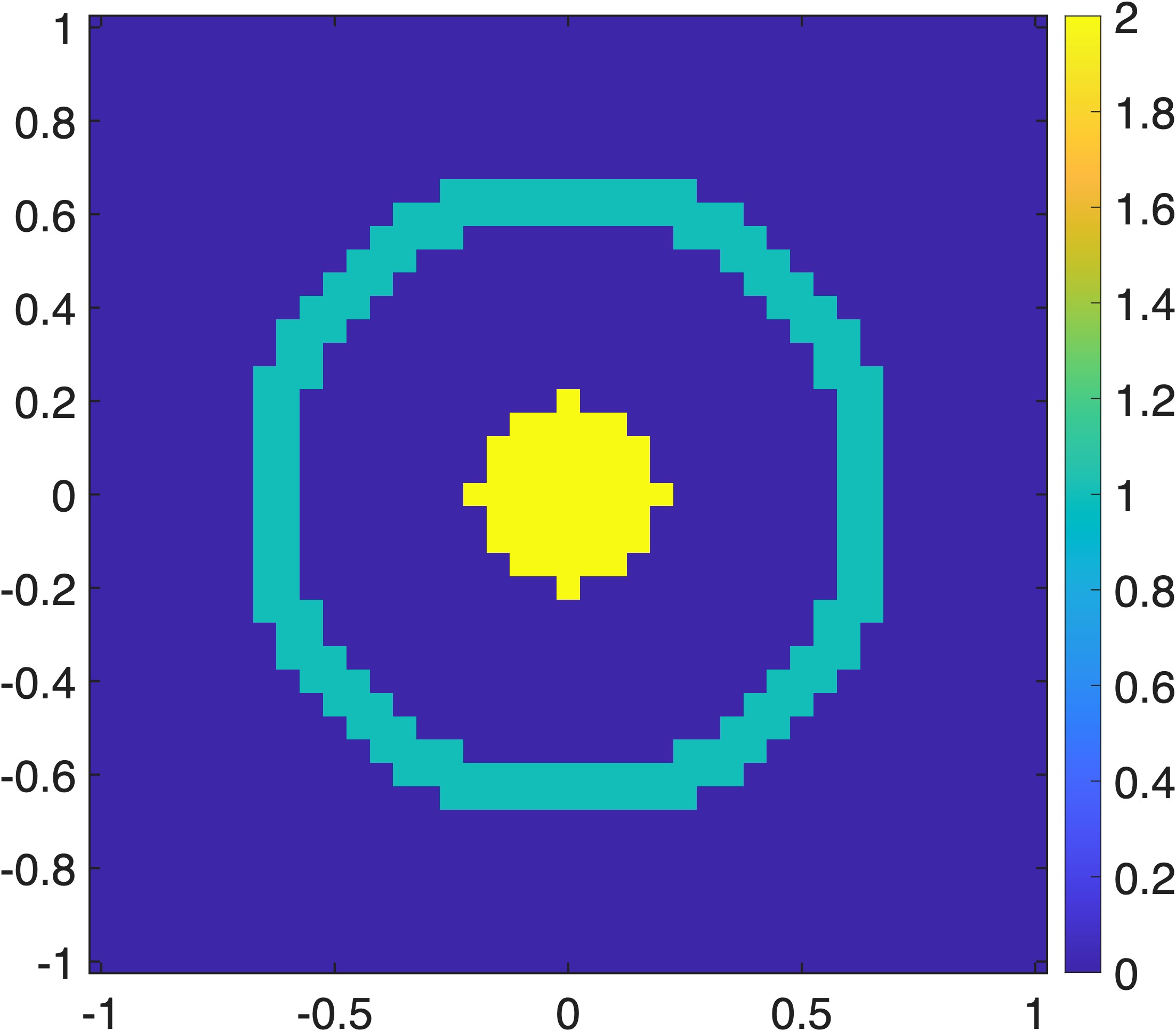}}
	\quad
	\subfloat[\label{fig3b}Reconstructed first velocity component $u^0_{\mathrm{comp},1}$]{\includegraphics[width=.25\textwidth]{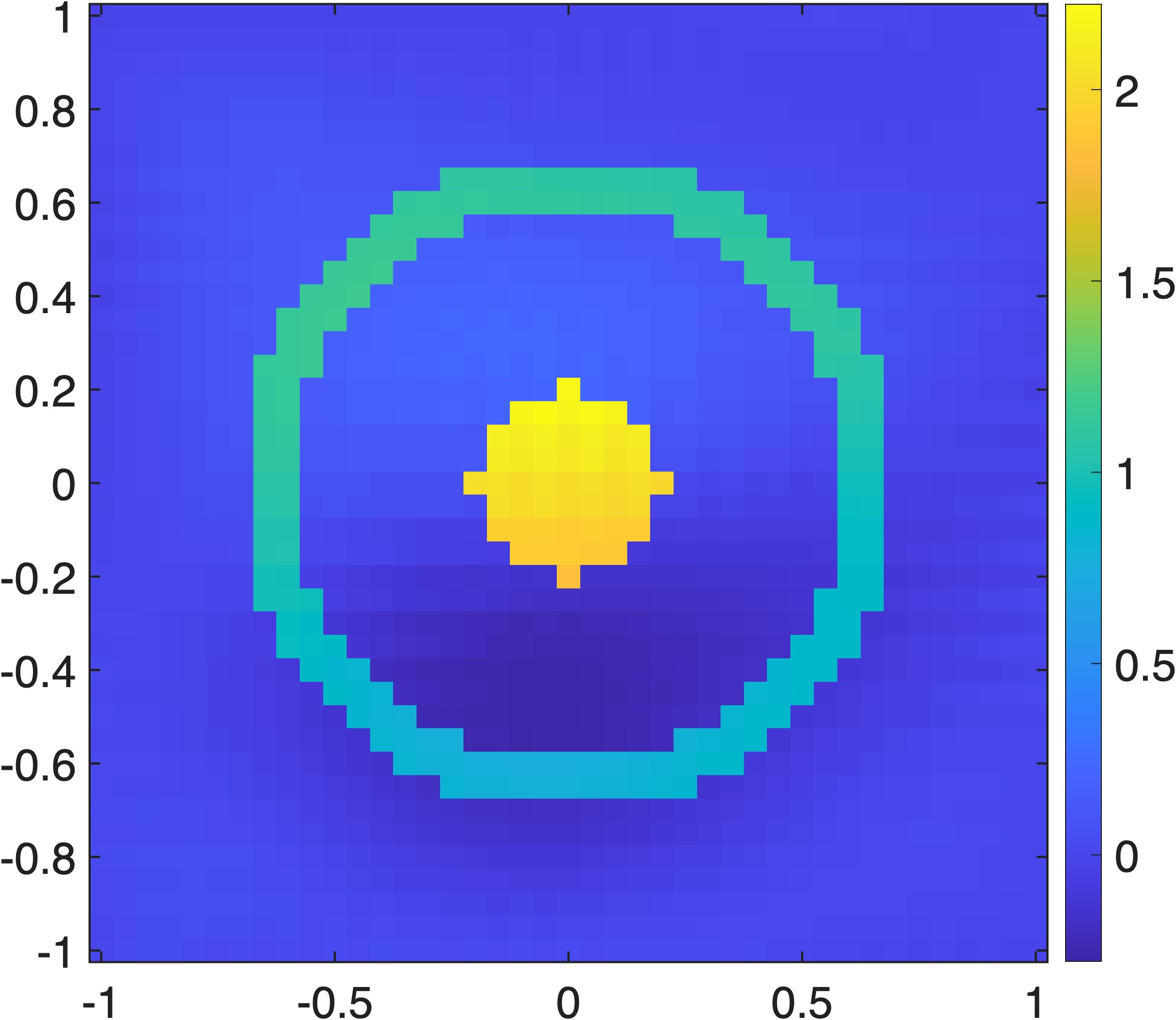}}
	\quad
	\subfloat[\label{fig3c}Relative error of the first component $\frac{|u^0_{{\rm true},1} - u^0_{{\rm comp},1}|}{\|u^0_{{\rm true},1}\|_{L^\infty(\Omega)}}$
]{\includegraphics[width=.25\textwidth]{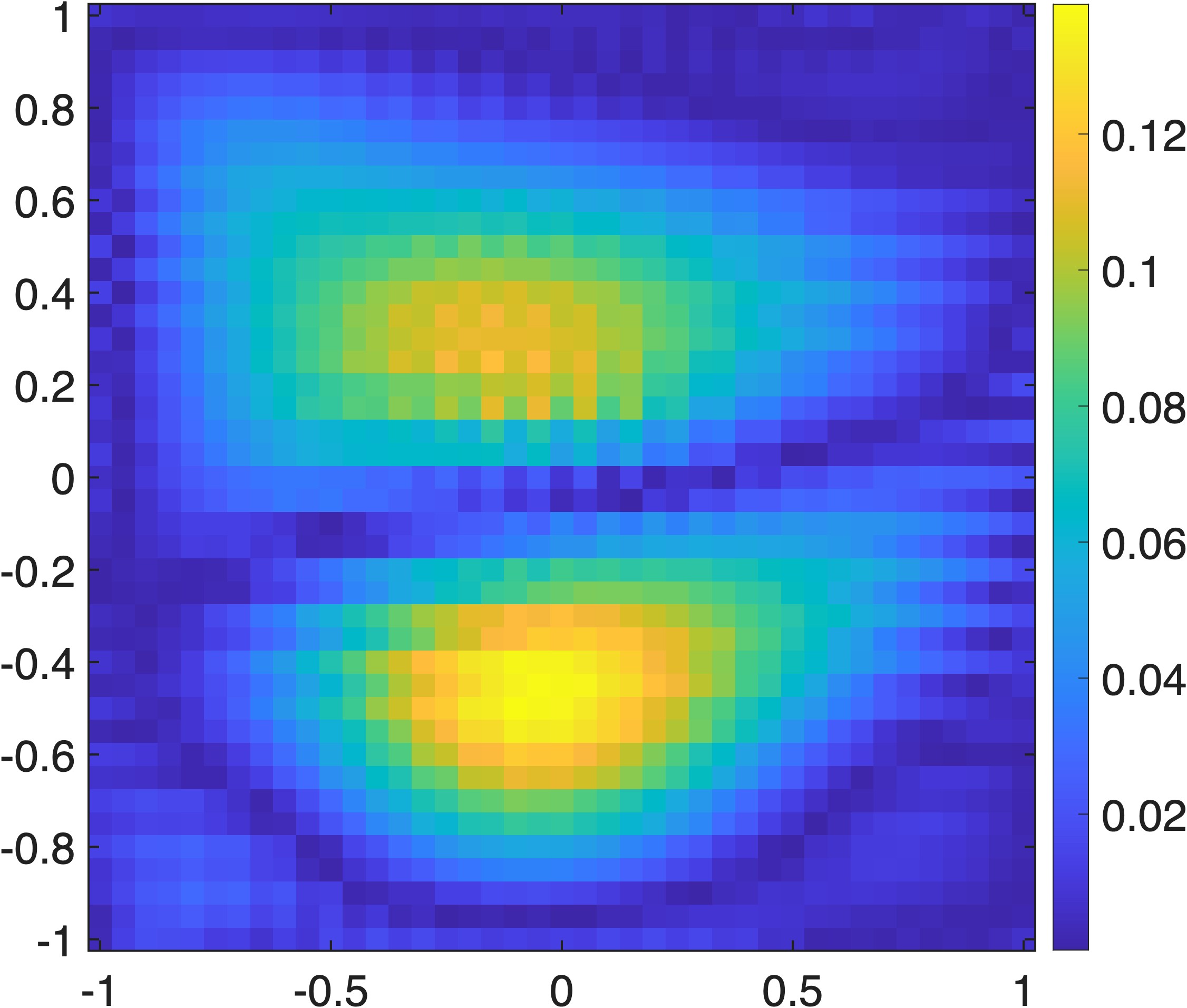}}
	
	\subfloat[\label{fig3d} True second velocity component $u^0_{\mathrm{true},2}$]{\includegraphics[width=.25\textwidth]{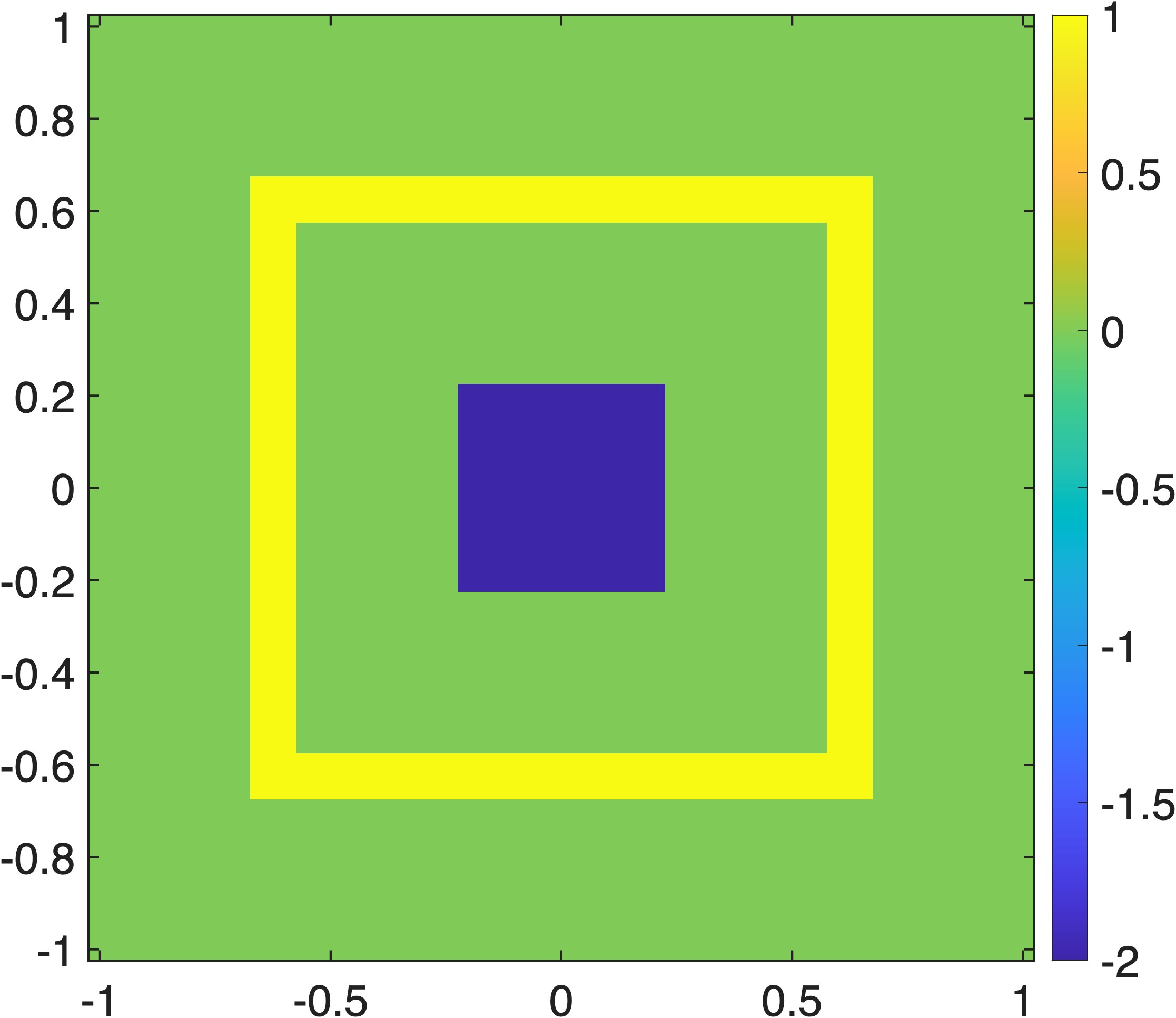}}
	\quad
	\subfloat[\label{fig3e}Reconstructed second velocity component $u^0_{\mathrm{comp},2}$]{\includegraphics[width=.25\textwidth]{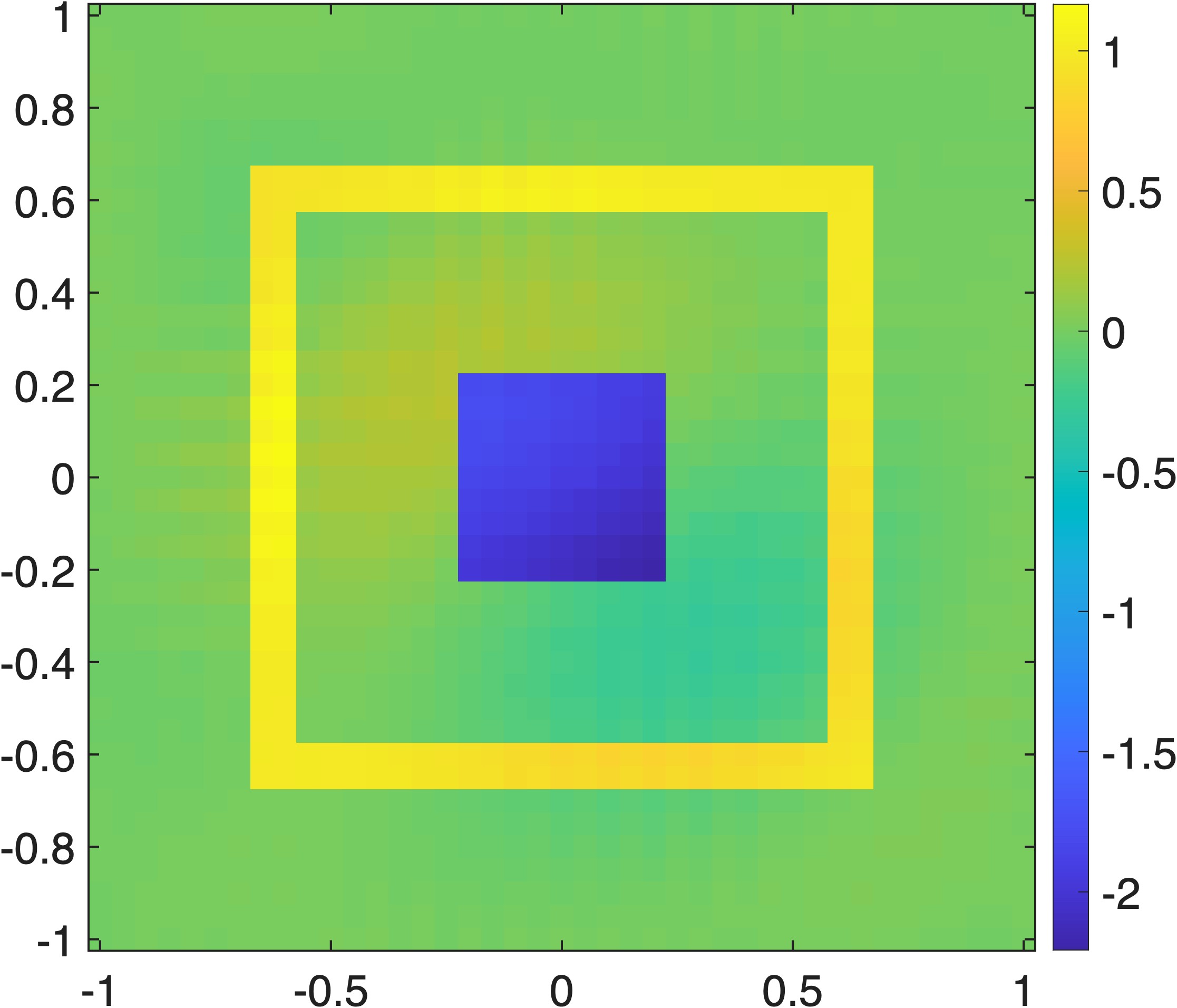}}
	\quad
	\subfloat[\label{fig3f}Relative error of the second component $\frac{|u^0_{{\rm true},2} - u^0_{{\rm comp},2}|}{\|u^0_{{\rm true},2}\|_{L^\infty(\Omega)}}$]{\includegraphics[width=.25\textwidth]{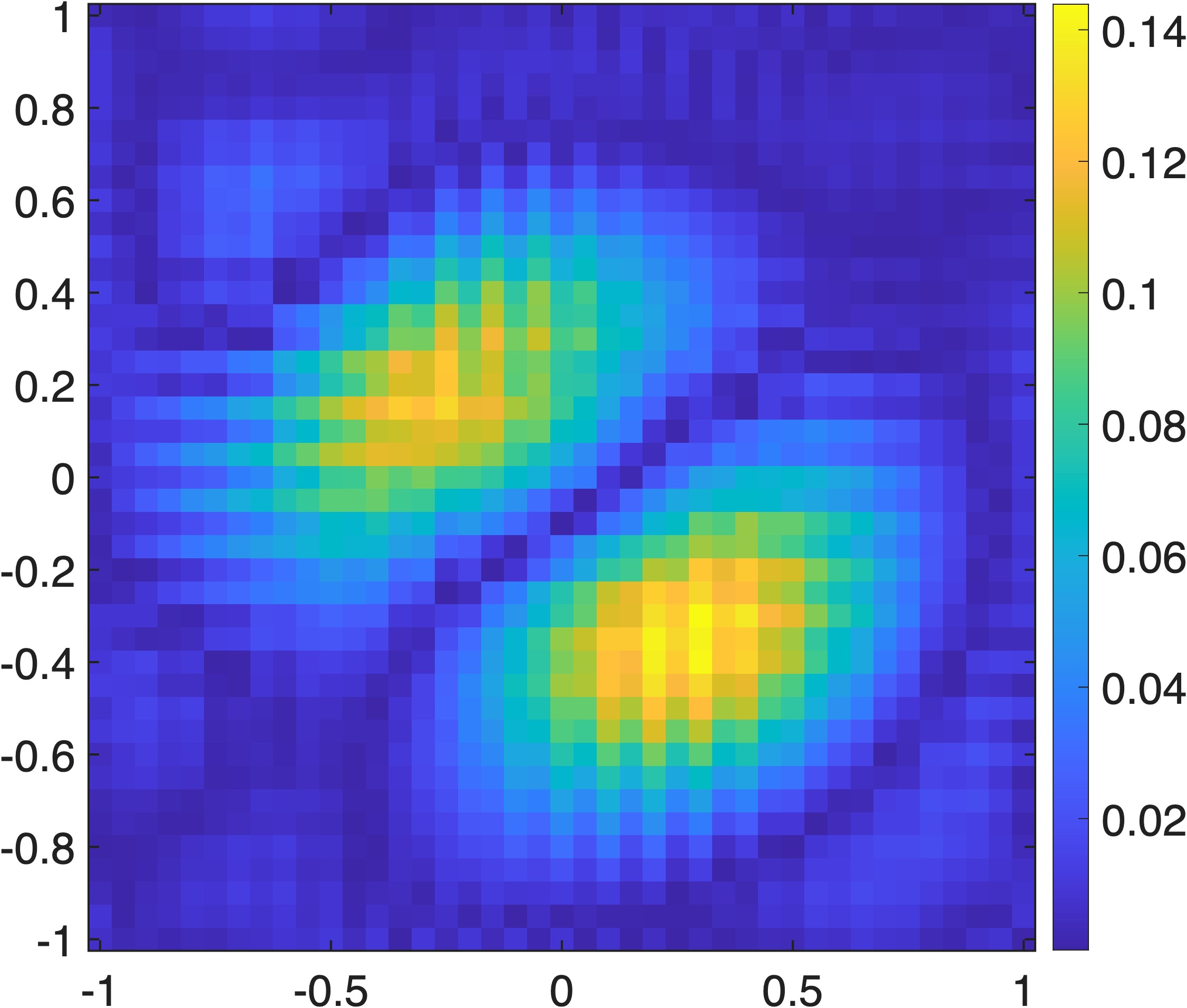}}
	
	\subfloat[\label{fig3g} True initial velocity field $\mathbf{u}^0_{\mathrm{true}}$]{\includegraphics[width=.25\textwidth]{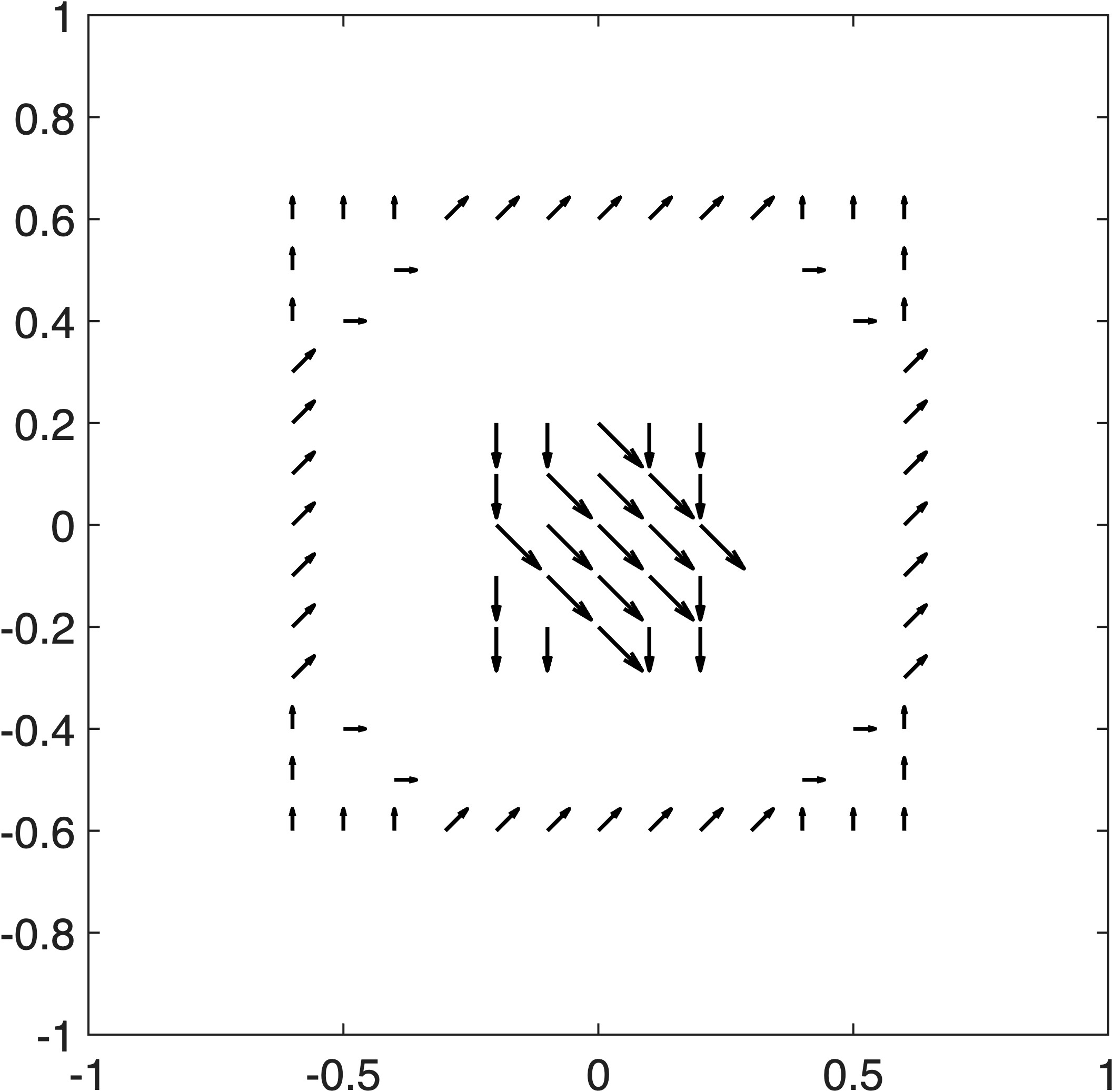}}
	\quad
	\subfloat[\label{fig3h} Reconstructed initial velocity field $\mathbf{u}^{0}_{\mathrm{comp}}$]{\includegraphics[width=.25\textwidth]{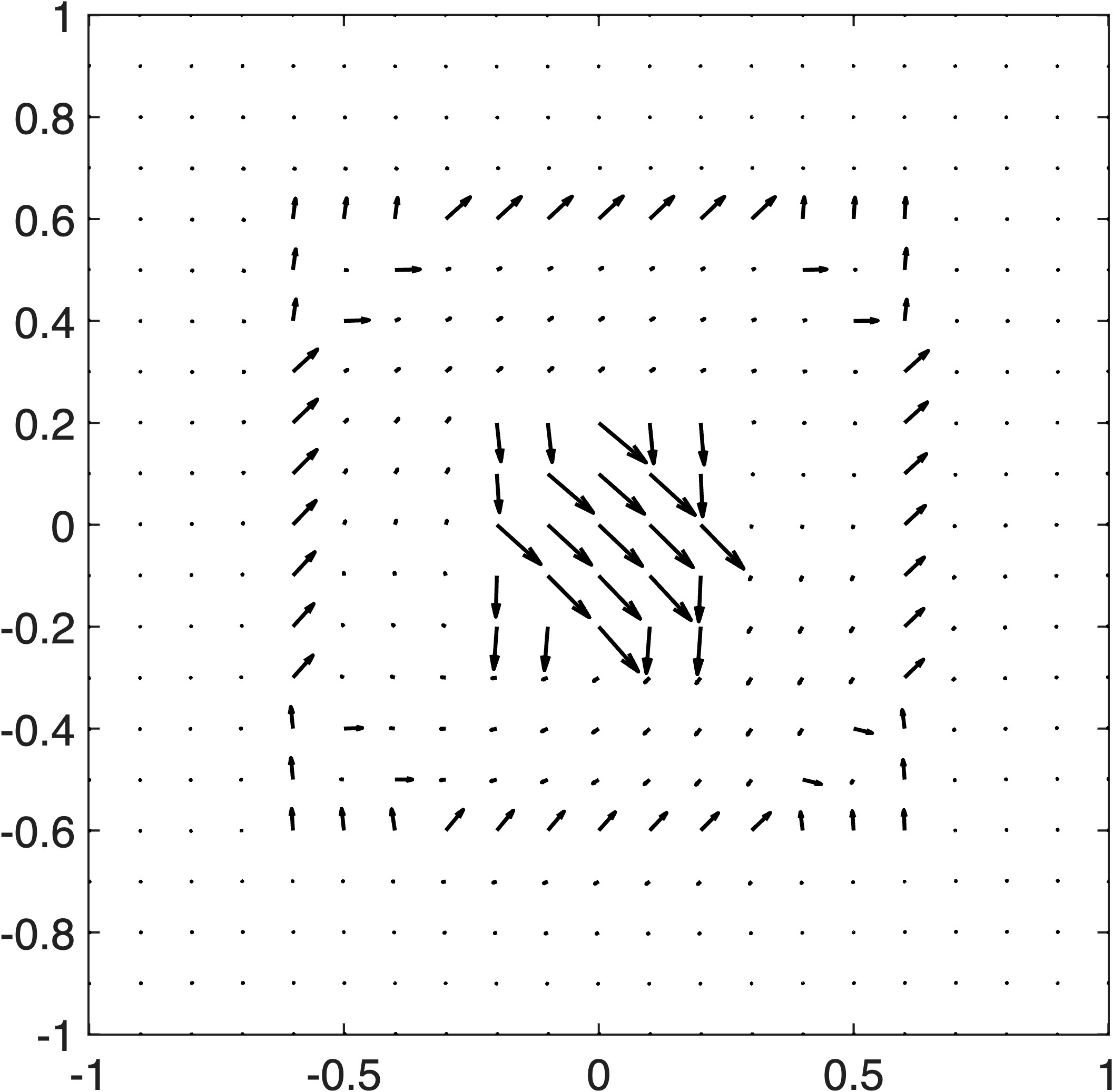}}
	\quad
	\subfloat[\label{fig3i} Residual of \eqref{time_red} for $\bU^{(k)}$ generated by Algorithm \ref{alg}]{\includegraphics[width=.25\textwidth]{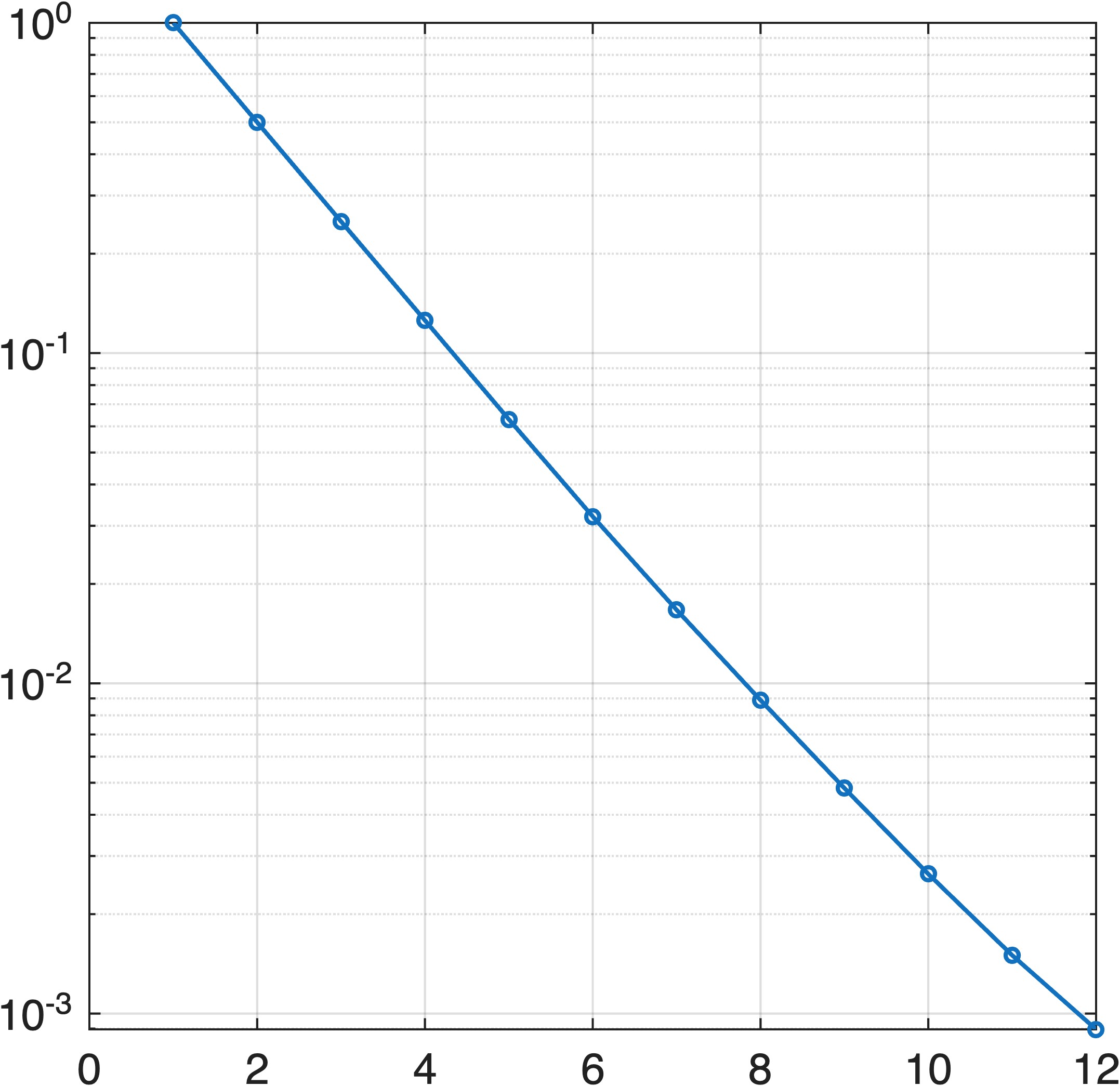}}

	\caption{\label{test3} Test 3: true and reconstructed initial velocity fields. Figures (a)--(c) show the true first component, its reconstruction, and the corresponding relative error, respectively. Figures (d)--(f) show the true second component, its reconstruction, and the corresponding relative error, respectively. Figures (g)--(h) display the true and reconstructed initial velocity fields. Figure (i) shows the residual of \eqref{time_red} for $\bU^{(k)}$ generated by Algorithm~\ref{alg}, which decreases steadily throughout the iteration.}
\end{figure}

The results shown in Figure~\ref{test3} demonstrate that the proposed method remains effective for more intricate initial data containing multiple layers, disconnected active regions, and sign changes. Figure~\ref{fig3a} and Figure~\ref{fig3b} show that the first component is reconstructed with the correct radial organization: both the outer annular structure and the inner disk are clearly identified, with their locations and relative magnitudes reasonably preserved. Likewise, Figure~\ref{fig3d} and Figure~\ref{fig3e} indicate that the second component successfully recovers the square outer ring together with the negative inner core, which confirms that the method is capable of distinguishing structures of different geometry and opposite sign within the same test. Although the reconstructed profiles exhibit noticeable smoothing and mild background artifacts, the essential topology and contrast pattern of the true solution are retained. This is further supported by the relative error distributions in Figure~\ref{fig3c} and Figure~\ref{fig3f}, where the displayed maxima are approximately $0.13$ and $0.14$, respectively. In addition, the comparison between Figure~\ref{fig3g} and Figure~\ref{fig3h} shows that the reconstructed initial velocity field captures the main directional organization of the true field. Finally, Figure~\ref{fig3i} shows a steady decay of the residual of \eqref{time_red} for $\bU^{(k)}$, providing further numerical evidence for the stability and robustness of the damped Picard iteration. Overall, this example illustrates that the proposed method can successfully recover composite structures with multiple scales, sharp interfaces, and sign-changing behavior from noisy boundary observations.

\section{On the stability of the time reduction and the choice of time basis}\label{sec:discuss}

We detail how temporal projection and truncation act as a low-pass filter, justify the choice of the exponentially weighted Legendre basis, and discuss a scaled weight for large $T$.

\subsection{Why time-dimensional reduction mitigates ill-posedness}
It is well known that the inverse problem under consideration is ill-posed: small perturbations in the data can induce arbitrarily large perturbations in the reconstruction. We now explain how the time-dimensional reduction helps to reduce the degree of ill-posedness in practice, as demonstrated by the acceptable numerical examples above.

Recall that in the time-dimensional reduction framework, we only need the ``low oscillation" time-projected boundary data
\[
    \mathbf{f}_m(\mathbf{x})
    \;=\;
    \partial_\nu \mathbf{u}_m(\mathbf{x})
    \;=\;
    \int_{0}^{T} e^{-2t}\,\mathbf{f}(\mathbf{x},t)\,\Psi_m(t)\,dt,
    \qquad m=0,1,\dots,N,
\]
while the higher oscillation components when $m > N$ do not contribute to the computational procedure.
Hence, when the measured data $\mathbf{f}$ contains noise, this projection-truncation step acts as a \emph{temporal low-pass filter}. In other words, we only feed into the inversion those time components of the data that survive the weighted projection onto the first $N$ modes. This averaging and truncation reduce the effective noise level in the data-misfit and, hence, relax the practical severity of the ill-posedness.

\subsection{Why use the exponentially weighted time basis} \label{subsec:exponent}
In this paper, we set $\Psi_n(t)=e^{t}Q_n(t)$, where $\{Q_n\}_{n\ge 0}$ are (shifted) Legendre polynomials. Although the factor $e^{t}$ ``cancels" in the $L^2_{e^{-2t}}(0,T)$ inner product, it is essential when time derivatives appear. Indeed,
\[
\Psi_n'(t)=\frac{d}{dt}\big(e^{t}Q_n(t)\big)=e^{t}\big(Q_n(t)+Q_n'(t)\big),
\]
which is not identically zero on $(0,T)$. Consequently, in expressions such as \eqref{4.3}, every coefficient $\bu_n$ contributes through the term $\sum_{n=0}^{N} \bu_n(\mathbf{x})\,\Psi_n'(t)$, and this information propagates into the time-reduced model \eqref{4.5}.

\begin{figure}[h!]
\centering
	\subfloat[True first velocity component $u^0_{\mathrm{true},1}$]{\includegraphics[width=.25\textwidth]{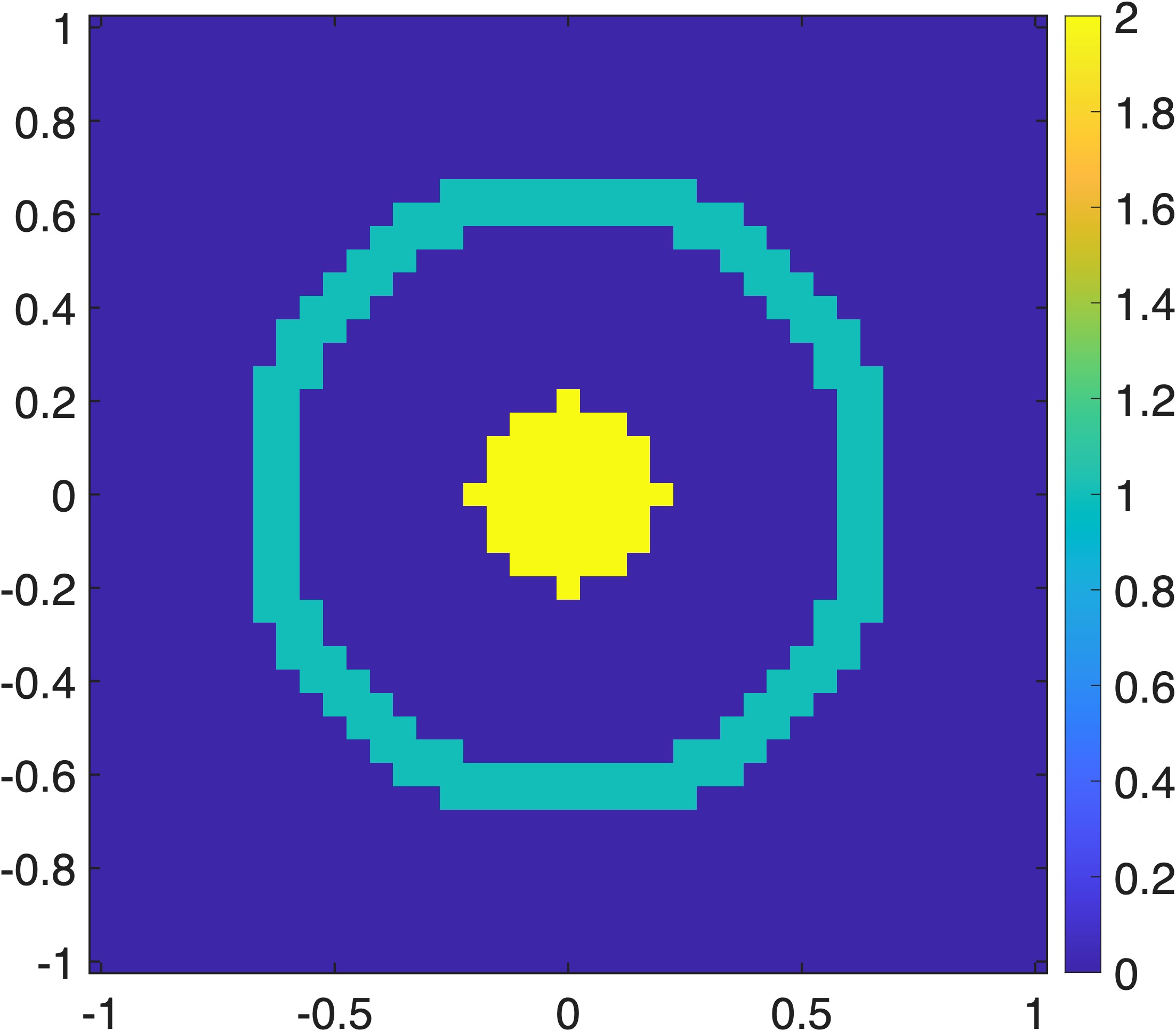}}
	\quad
	\subfloat[\label{fig4b}Reconstructed first velocity component $u^0_{\mathrm{comp},1}$]{\includegraphics[width=.25\textwidth]{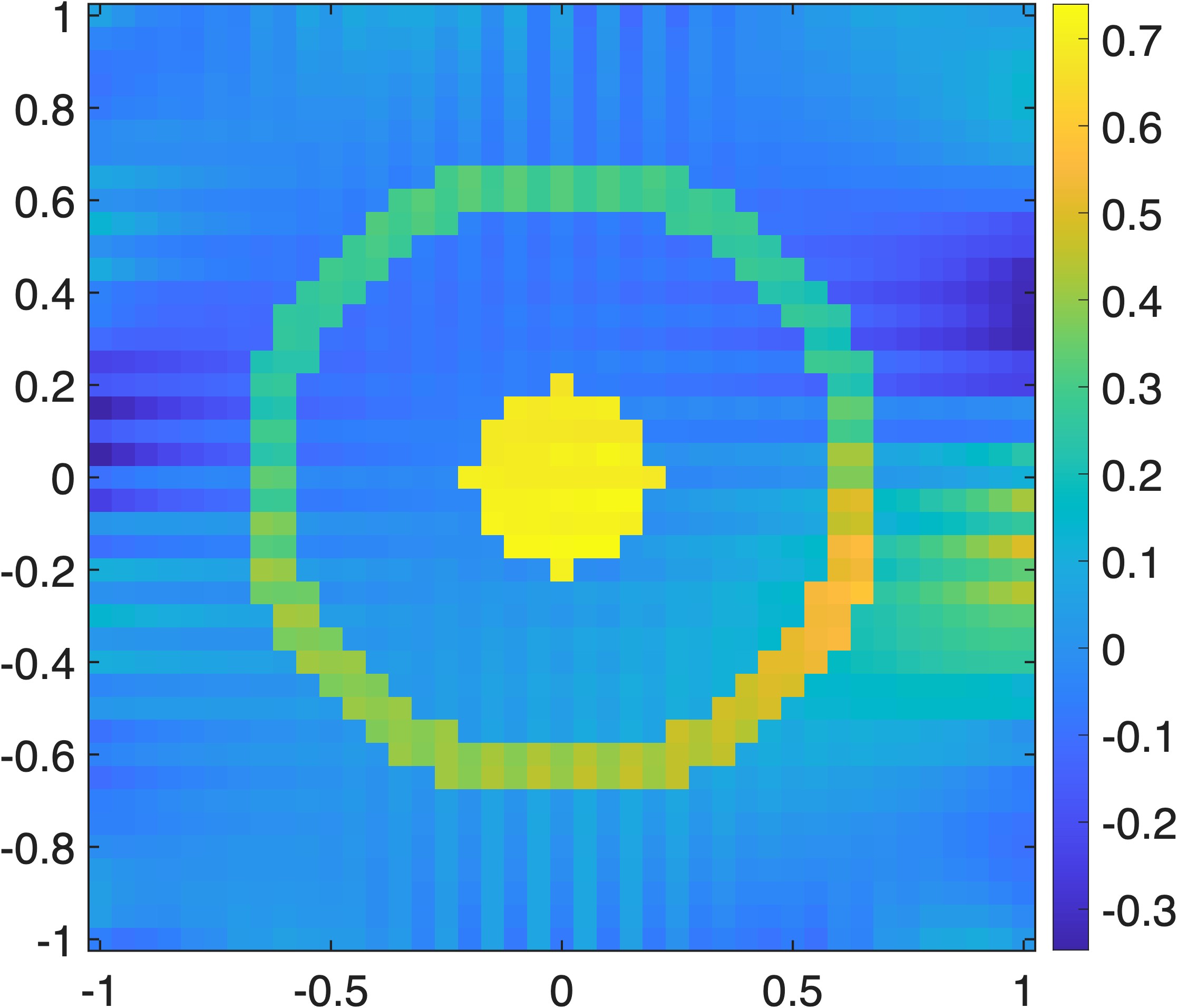}}
	\quad
	\subfloat[\label{fig4c}Relative error of the first component $\frac{|u^0_{{\rm true},1} - u^0_{{\rm comp},1}|}{\|u^0_{{\rm true},1}\|_{L^\infty(\Omega)}}$
]{\includegraphics[width=.25\textwidth]{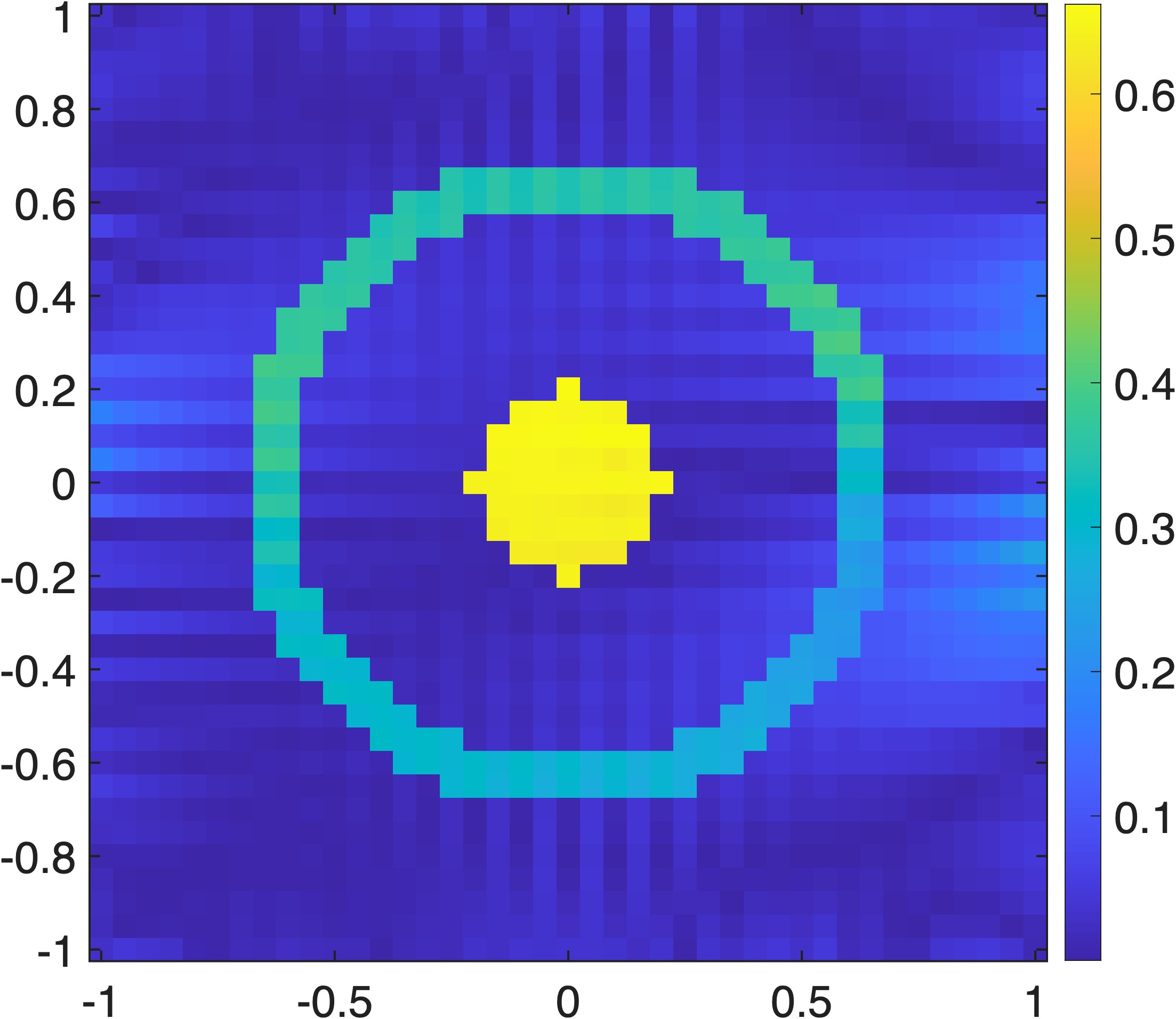}}
	
	\subfloat[ True second velocity component $u^0_{\mathrm{true},2}$]{\includegraphics[width=.25\textwidth]{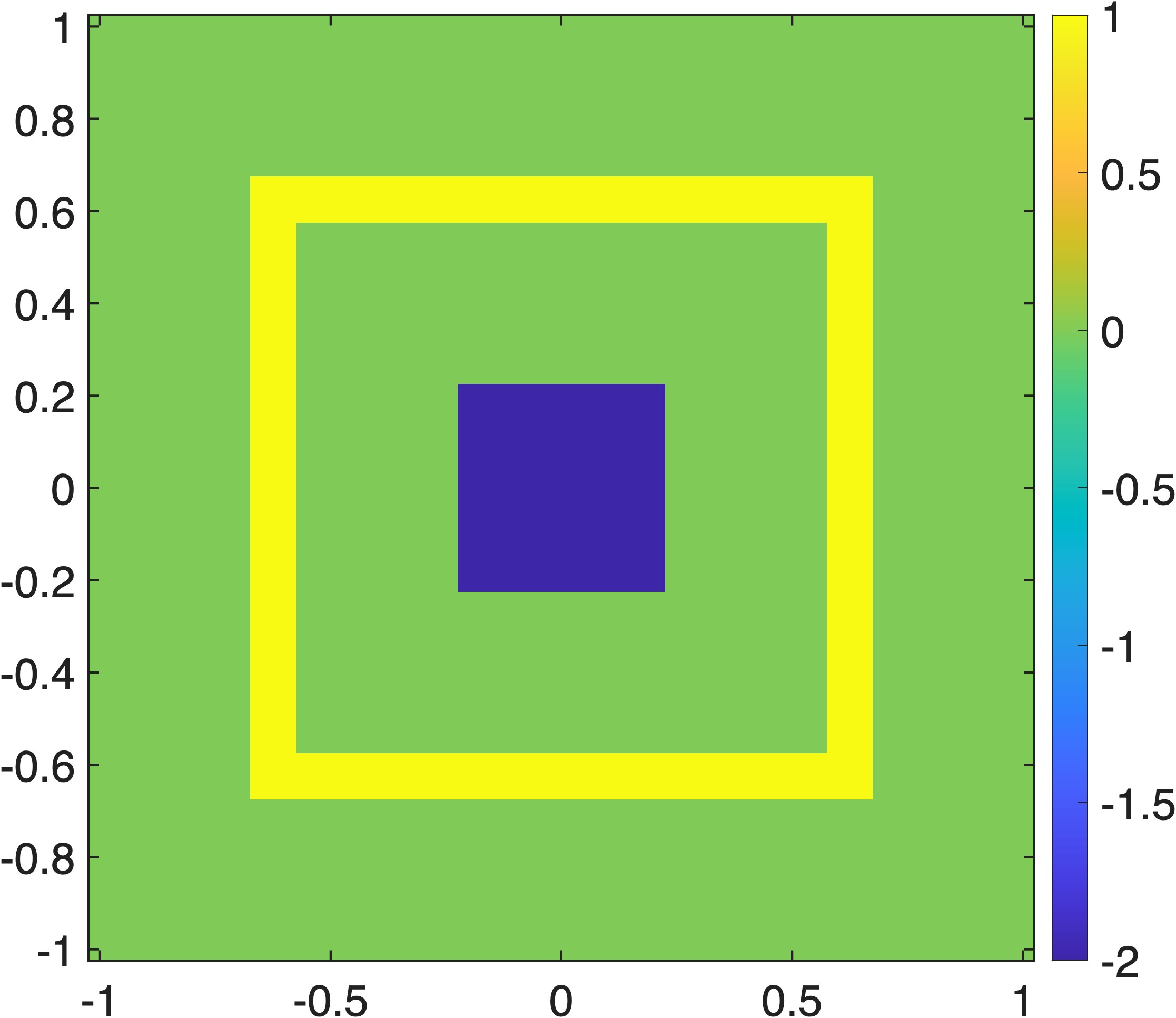}}
	\quad
	\subfloat[\label{fig4e}Reconstructed second velocity component $u^0_{\mathrm{comp},2}$]{\includegraphics[width=.25\textwidth]{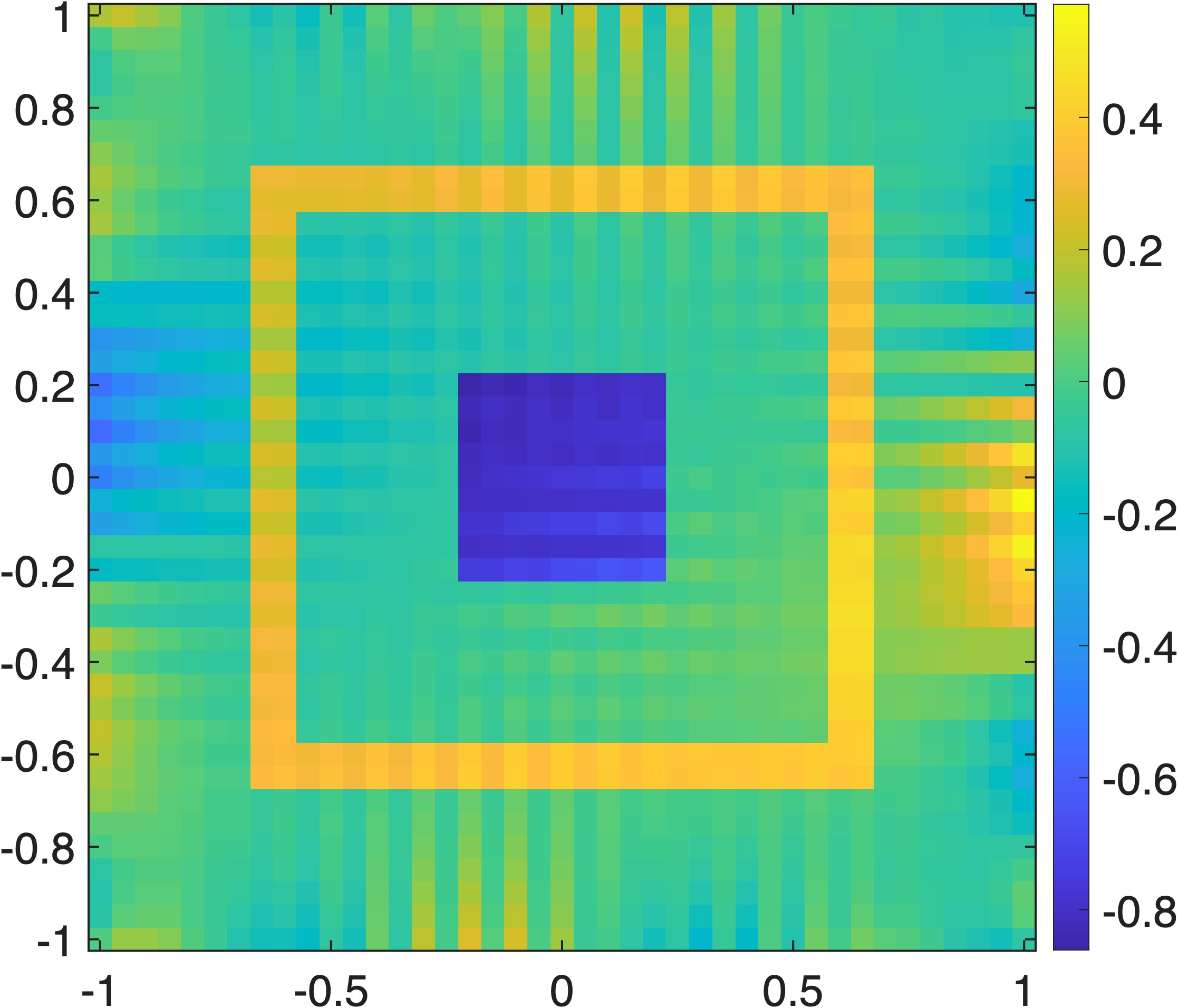}}
	\quad
	\subfloat[\label{fig4f}Relative error of the second component $\frac{|u^0_{{\rm true},2} - u^0_{{\rm comp},2}|}{\|u^0_{{\rm true},2}\|_{L^\infty(\Omega)}}$]{\includegraphics[width=.25\textwidth]{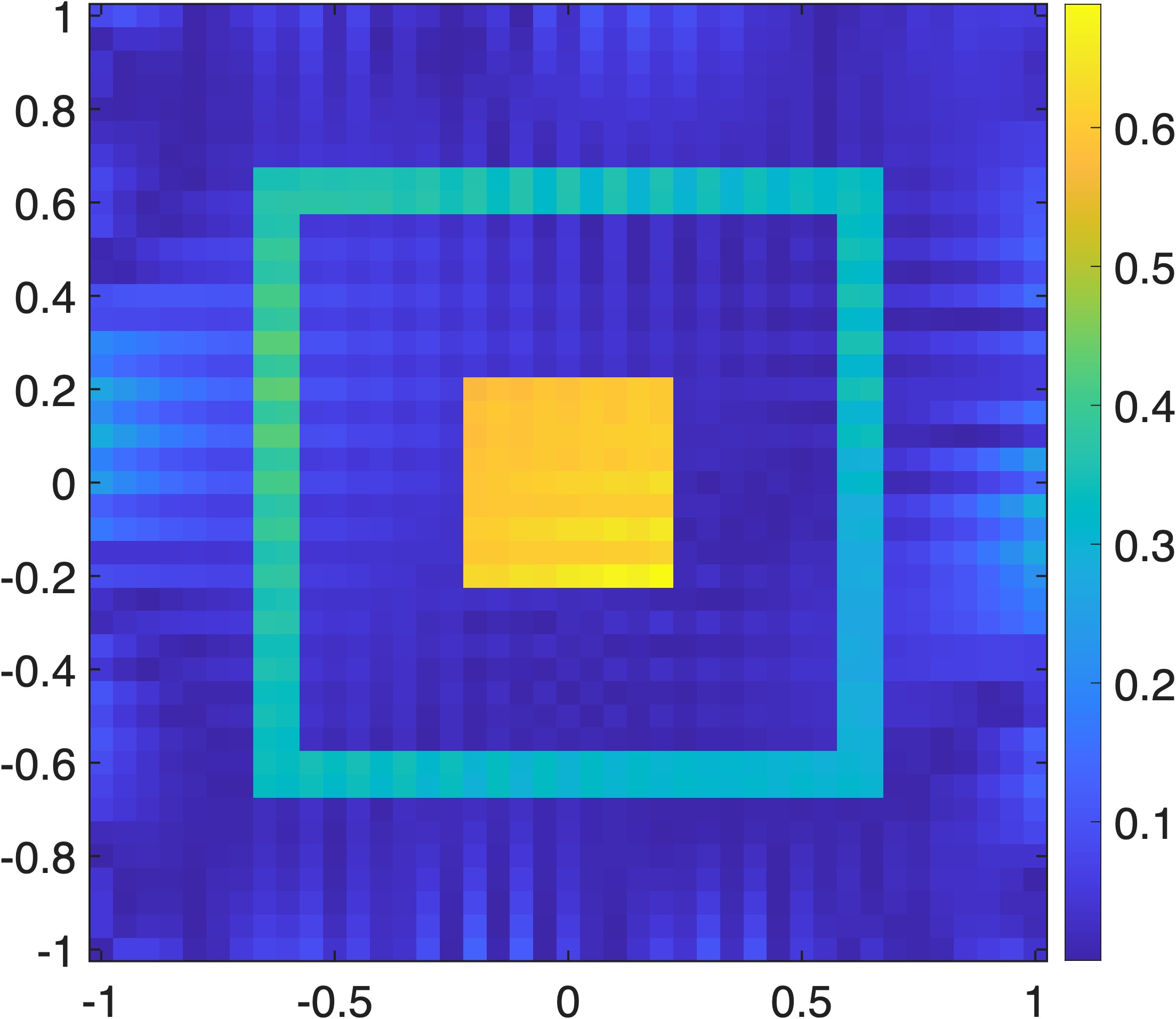}}
	\caption{\label{test3_noWeight} Test 3 using the unweighted Legendre basis $\{Q_n\}_{n\ge 0}$ in $L^2(0,T)$. Figures (a)--(c) show the true first component, its reconstruction, and the corresponding relative error, respectively. Figures (d)--(f) show the true second component, its reconstruction, and the corresponding relative error, respectively.}
\end{figure}

By contrast, using the unweighted Legendre basis $\{Q_n\}_{n\ge 0}$ gives $Q_0'(t)=0$ for all $t$, as $Q_0$ is a constant. The corresponding coefficient $\bu_0$ can be underrepresented in the derivative equations, which may weaken coupling between modes and degrade numerical conditioning. The exponential weighting avoids this loss by ensuring that even the lowest time mode remains visible to the derivative terms.
This observation is confirmed by recomputing Test~3 with the time-dimensional reduction method using the unweighted basis $\{Q_n\}_{n\ge 0}$ in $L^2(0,T)$ in place of $\{\Psi_n\}_{n\ge 0}$; as shown in Figure~\ref{test3_noWeight}, the exponentially weighted basis $\{\Psi_n\}_{n\ge 0}$ yields visibly better reconstructions.

Figure~\ref{test3_noWeight} illustrates the reconstruction in Test~3 when the unweighted Legendre basis $\{Q_n\}_{n\ge0}$ in $L^2(0,T)$ is used in place of the exponentially weighted basis $\{\Psi_n\}_{n\ge0}$. Figure~\ref{fig4b} and Figure~\ref{fig4e} show that, although the main annular and square-ring structures are still identifiable, the reconstructed amplitudes are significantly damped and the interfaces are much more diffuse than those obtained with the weighted basis in Figure~\ref{test3}. In particular, the central disk and the inner square are recovered only in a blurred manner, and noticeable background contamination appears throughout the domain. This degradation is reflected in the relative error distributions shown in Figure~\ref{fig4c} and Figure~\ref{fig4f}, where the displayed maxima are approximately $0.65$ and $0.67$, respectively, which are substantially larger than those observed in Figure~\ref{test3}. Compared with the results produced by the exponentially weighted basis, the unweighted basis yields weaker contrast, less accurate amplitudes, and poorer resolution of sharp interfaces. These numerical results therefore support the use of the exponentially weighted basis $\{\Psi_n\}_{n\ge0}$ in the time-dimensional reduction procedure.

\subsection{On the exponential weight for large $T$}
A natural concern is that the factor $e^{t}$ can become very large when $T$ is large. To avoid this, we may replace $e^{t}$ by a scaled weight $e^{\alpha t}$ with $\alpha=1/T$ (or more generally $\alpha=c/T$ for a modest constant $c$) and work in $L^2_{e^{-2\alpha t}}(0,T)$. In this case $\max_{t\in[0,T]} e^{\alpha t}=e^{\alpha T}=e^{c}$, so the weight remains uniformly bounded as $T$ increases. The construction of the time-dimensional reduction and all subsequent steps in this paper remain the same under this scaling; only the value of $\alpha$ is adjusted to control the dynamic range.

\section{Concluding Remarks}\label{sec:conclusion}

This work presents a new computational approach to a targeted inverse initial data problem for the compressible anisotropic Navier--Stokes system. The objective is to reconstruct the initial velocity field $\bu^0$ from noisy boundary observations, while the density $\rho(\x,t)$, the pressure $p(\x,t)$, the anisotropic viscosity tensor $\bmu$, and the body force $\bF(\x,t)$ are assumed known. The proposed approach is based on an exponentially weighted Legendre time-dimensional reduction, which transforms the original time-dependent inverse problem into a coupled system of time-independent elliptic equations for the Fourier coefficients of the velocity field. This reduction makes the problem substantially more tractable both analytically and computationally.

To solve the resulting time-reduced model, we combine quasi-reversibility with a damped Picard iteration. The numerical results demonstrate that the proposed method yields accurate and stable reconstructions across several test cases, including examples with localized structures, non-axis-aligned features, multiple disconnected regions, and sign changes. The experiments also show robustness with respect to noisy boundary data, which indicates that the method has promising potential for inverse problems in realistic settings where the initial state is not directly observable.

An important open question is whether the sequence $\bigl\{\mathbf{U}^{(k)}\bigr\}_{k\geq 0}$ defined in Section~\ref{sec:Picard} converges to the true solution of the time-reduced model~\eqref{time_red}. A positive answer would constitute a significant advance in the analysis of the proposed method. One possible route toward such a result would be to incorporate Carleman weights into the functional $J_{{\bf U}^{(k)}}^{\epsilon,N}$ and to combine this construction with suitable Carleman estimates. However, such estimates are not currently available in the fully anisotropic setting considered here.

If one restricts instead to the isotropic case, in particular to the Lam\'e system, then existing Carleman estimates become applicable; see, e.g., \cite{ImanuvilovLorenziYamamoto2022, Isakov2007}. In that setting, several recently developed techniques, including Carleman convexification \cite{VoKlibanovNguyen:IP2020, KlibanovIoussoupova:SMA1995, KlibanovNguyenTran:JCP2022, LeLeNguyen:2024}, the Carleman contraction principle \cite{LeCON2023, LeNguyen:jiip2022, Nguyen:AVM2023, NguyenNguyenVu2026}, and the Carleman--Newton method \cite{AbhishekLeNguyenKhan, LeNguyenTran:CAMWA2022}, may be adapted to construct a globally convergent sequence $\bigl\{\mathbf{U}^{(k)}\bigr\}_{k\geq 0}$ that approximates the true solution of the time-reduced model~\eqref{time_red}. In such a case, both analytical and numerical global solvability of Problem~\ref{isp} may become achievable. Extending these ideas to the anisotropic setting remains a central direction for future research.


\begin{thebibliography}{10}
\bibitem{Liu2024}
G.~Liu.
\newblock Determining the viscosity function from the boundary measurements for
  the stokes and the navier--stokes equations.
\newblock {\em Inverse Problems}, 40(12):125011, dec 2024.
\end{thebibliography}

\end{document}